\documentclass[letterpaper]{article}
\usepackage{aaai19}
\usepackage{times}
\usepackage{helvet}
\usepackage{courier}
\usepackage{amsmath}
\usepackage{amssymb}
\usepackage{algorithm}
\usepackage{graphicx}
\usepackage{subfigure}
\usepackage{algorithmic}
 \usepackage{natbib}
\newcommand{\EE}{\mathbb{E}}

\newcommand{\argmin}{\mathop{\mathrm{argmin}}}
\usepackage{theorem}
\newtheorem{lemma}{Lemma}

\newtheorem{theorem}{Theorem}

\newtheorem{corollary}{Corollary}

\begin{document}
%
\title{Non-ergodic Convergence Analysis of Heavy-Ball Algorithms}
\author{Tao Sun$^{\dagger}$, Penghang Yin$^{\ddagger}$, Dongsheng Li$^{\sharp}$,
Chun Huang$^{\sharp}$,
Lei Guan$^{\sharp}$,
Hao Jiang$^{\sharp}$\\
$^{\dagger}$Department of Mathematics, National University of Defense Technology, Changsha, Hunan, China. \\
$^{\ddagger}$Department of Mathematics, University of California, Los Angeles, USA. \\
$^{\sharp}$College  of Computer, National University of Defense Technology, Changsha, Hunan, China.\\
Emails: \textit{nudtsuntao@163.com}, \textit{yph@ucla.edu}, \textit{\{dsli,chunhuang\}@nudt.edu.cn}\\ \textit{guanleics@gmail.com}, \textit{haojiang@nudt.edu.cn}
}
\maketitle
\begin{abstract}
\begin{quote}
In this paper, we revisit the convergence of the Heavy-ball method, and present improved convergence complexity results in the convex setting. We provide the first non-ergodic $O(1/k)$ rate result of the Heavy-ball algorithm with constant step size for coercive objective functions. For objective functions satisfying a relaxed strongly convex condition, the linear convergence is established under weaker assumptions on the step size and inertial parameter than made in the existing literature.  We extend our results to multi-block version of the algorithm with both the cyclic and stochastic update rules. In addition, our results can also be extended to decentralized optimization, where the ergodic analysis is not applicable.
\end{quote}
\end{abstract}

\section{Introduction}
In this paper, we study the Heavy-ball algorithm first proposed by \cite{polyak1964some}, for solving the following unconstrained minimization problem
\begin{equation}\label{model}
    \min_{x\in \mathbb{R}^m} f(x),
\end{equation}
where $f$ is convex and differentiable, and $\nabla f$ is Lipschitz continuous with constant $L$. The Heavy-ball method iterates
 \begin{align}\label{scheme}
     x^{k+1}=x^k-\gamma_k \nabla f(x^k)+\beta_k(x^k-x^{k-1}),
\end{align}
where $\gamma_k$ is the step size and $\beta_k$ is the inertial parameter.
Different from the gradient descent algorithm, the sequence generated by Heavy-ball method is not Fej\'{e}r monotone due to the inertial term $\beta_k(x^k-x^{k-1})$.  This poses a challenge in proving the convergence rate of $\{f(x^k)\}_{k\geq 0}$ in the convex case.  In the existing literature, the sublinear convergence rate of the Heavy-ball has been proved only in the sense of ergodicity.

When the objective function is twice continuously differentiable and strongly convex (i.e., almost quadratic), the Heavy-ball method provably converges linearly. Under a weaker assumption that the objective function is nonconvex but Lipschitz differentiable, \cite{zavriev1993heavy} proved   that the sequence generated by the Heavy-ball method will converge to a critical point, yet without specifying the convergence rate. The smoothness of objective function is crucial for convergence of the Heavy-ball. Indeed, it can be divergent for a strongly convex but nonsmooth function as suggested by \cite{lessard2016analysis}.   Different from the classical gradient descent methods, the Heavy-ball algorithm fails to generate a Fej\'{e}r monotone sequence. In the convex and smooth case, the only result about convergence rate, to our knowledge, is the ergodic $O(1/k)$ rate in terms of the objective value \citep{ghadimi2015global}, i.e., $f\left(\frac{\sum_{i=1}^k x^i}{k}\right)-\min f\sim O\left(\frac{1}{k}\right)$. The linear convergence of Heavy-ball algorithm was proved under the strongly convexity assumption by \cite{ghadimi2015global}. But the authors imposed a restrictive assumption on the inertial parameter $\beta_k$. Specifically, when the strongly convex constant is tiny, the convergence result holds only for a small range of $\beta_k$ values. By incorporating the idea of proximal mapping, the inertial proximal gradient algorithm (iPiano) was proposed in \citep{ochs2014ipiano}, whose convergence in nonconvex case was thoroughly discussed. Locally linear convergence of iPiano and Heavy-ball method was later proved in \citep{ochs2016local}. In the strongly convex case, the linear convergence was proved for iPiano  with fixed $\beta_k$ \citep{ochs2015ipiasco}. In the paper \citep{pock2016inertial}, inertial Proximal Alternating Linearized Minimization (iPALM) was introduced as a variant of iPiano for solving the two-block regularized problem. \cite{xu2013block} analyzed the Heavy-ball algorithm in tensor minimization problems. Stochastic versions of heavy-ball have also been introduced \citep{loizou2017momentum, loizou2017linearly}. A multi-step heavy-ball algorithm was analyzed in \citep{liang2016multi}. The inertial methods are also developed and studied in the operator research by \cite{Combettes2017Quasinonexpansive}.  None of the aforementioned Heavy-ball based algorithms, however, provides a non-ergodic convergence rate.

\subsection{Contributions}
In this paper, we establish the \emph{first} non-ergodic $O(1/k)$ convergence result in general convex case. More precisely, we prove that $f(x^k)-\min f\sim O(\frac{1}{k})$ for convex  and coercive $f$ \footnote{We say $f$ is coercive, if $f(x)\rightarrow+\infty$ as $x\rightarrow+\infty$.}. Compared with existing result in \citep{ghadimi2015global}, ours allows a larger step size $\gamma_k$. We also prove a linear convergence result under a restricted strongly convex condition, weaker than the strong convexity assumption. In short, we make weaker assumptions on the step size, on the inertial parameter, as well as on the convexity of the objective function. The convergence of  multi-block extensions of Heavy-ball method is studied. The sublinear and linear convergence rates are proved for the cyclic and stochastic update rules, respectively. In addition, we extend our analysis to the decentralized  Heavy-ball method, where the ergodic analysis is not applicable. Our theoretical results are based on a novel Lyapunov function, which is motivated by  a modified dynamical system.


\subsection{A dynamical system interpretation}
It has been long known that the Heavy-ball method is equivalent to the  discretization of the following second-order ODE \citep{alvarez2000minimizing}:
\begin{equation}\label{ODE1}
    \ddot{x}(t)+\alpha \dot{x}(t)+\nabla f(x(t))=\textbf{0},~~t\geq 0,
\end{equation}
for some $\alpha>0$. In the case $\beta_k\equiv0$, the Heavy-ball method boils down to the standard gradient descent, which is known to be the discretization of the following first-order ODE
\begin{equation}\label{ODE2}
    \alpha \dot{x}(t)+\nabla f(x(t))=\textbf{0},~~t\geq 0.
\end{equation}
The dynamical system (\ref{ODE1}), however, misses essential information about relation between $\ddot{x}(t)$ and $\dot{x}(t)$. Specifically, if we replace $\ddot{x}(t)$ by $\frac{x^{k+1}-2x^{k}+x^{k-1}}{h^2}$ with $h$ being the discretization step size, then it holds that
\begin{align*}
    &\left\|\frac{x^{k+1}-2x^{k}+x^{k-1}}{h^2}\right\|\\
    &\quad\leq \frac{1}{h}\cdot\left( \left \|\frac{x^{k+1}-x^{k}}{h}\right \|+\left\|\frac{x^{k}-x^{k-1}}{h}\right\|\right).
\end{align*}
Since both $\frac{x^{k+1}-x^{k}}{h}$ and $\frac{x^{k}-x^{k-1}}{h}$  can be viewed as the discretization of $\dot{x}(t)$, we propose to modify (\ref{ODE1}) by adding the following constraint
 \begin{equation}\label{constraint}
    \|\ddot{x}(t)\|\leq \theta\|\dot{x}(t)\|,
 \end{equation}
 where $\theta>0$. In next section, we will devise a useful Lyapunov function by exploiting the additive constraint (\ref{constraint}) and establish the asymptotic non-ergodic sublinear convergence rate in the continuous setting. Finally, we will ``translate'' this analysis into that in discretized setting.

\section{Analysis of the dynamical system}
We analyze the modified dynamical system (\ref{ODE1}) + (\ref{constraint}). The existence of the solution is beyond the scope of this paper and will not be discussed here. Let us assume that  $f$ is coercive,  $\alpha>\theta$, and $f(x(0))-\min f>0$. We consider the Lyapunov function
\begin{align}\label{lyas}
\xi(t):= f(x(t))+\frac{1}{2}\|\dot{x}(t)\|^2-\min f \geq 0,
\end{align}
and refer the readers to the relevant equations (3)-(6). A direct calculation gives
\begin{align}\label{system-t1}
\dot{\xi}(t)&=\langle\nabla f(x(t)),\dot{x}(t)\rangle+\langle\ddot{x}(t),\dot{x}(t)\rangle\nonumber\\
&=-\alpha\|\dot{x}(t)\|^2,
\end{align}
which means $\{\xi(t)\}_{t\geq 0}$ is non-increasing. As a result,
 $$\sup_{t}\{f(x(t))-\min f\}\leq\sup_{t}\xi(t)\leq \xi(0).$$
By the coercivity of $f$, $\{x(t)\}_{t\geq 0}$ is bounded. Then by the continuity of $\nabla f$, $\{\nabla f(x(t))\}_{t\geq 0}$ is bounded; using \eqref{ODE1}, $\{\ddot{x}(t)+\alpha \dot{x}(t)\}_{t\geq 0}$ is also bounded. By the triangle inequality, we have
 \begin{align}
 \|\ddot{x}(t)+\alpha \dot{x}(t)\|&\geq \alpha\|\dot{x}(t)\|-\|\ddot{x}(t)\|\nonumber\\
 &\geq (\alpha-\theta)\|\dot{x}(t)\|.
 \end{align}
Since $\alpha>\theta$, we obtain the boundedness of $\{\dot{x}(t)\}_{t\geq 0}$; by \eqref{constraint}, $\{\ddot{x}(t)\}_{t\geq 0}$  is also bounded. Let $x^*\in \textrm{arg}\min f$, we have
 \begin{align}\label{rela}
&0\leq f(x(t))-f(x^*)\overset{a)}{\leq} \langle \nabla f(x(t)), x(t)-x^*\rangle\nonumber\\
&\quad\quad\overset{b)}{\leq}  \|\nabla f(x(t))\|\cdot\|x(t)-x^*\|\nonumber\\
&\quad\quad\overset{c)}{=} \|\ddot{x}(t)+\alpha \dot{x}(t)\|\cdot\|x(t)-x^*\|\nonumber\\
&\quad\quad\overset{d)}{\leq} (\|\ddot{x}(t)\|+\alpha \|\dot{x}(t)\|)\cdot\|x(t)-x^*\|\nonumber\\
&\quad\quad\overset{e)}{\leq}  (\alpha+\theta)\|\dot{x}(t)\|\cdot\|x(t)-x^*\|,
 \end{align}
 where $a)$ is due to the convexity of $f$; $b)$ is due to the Young's inequality; $c)$ is due to (\ref{ODE1}); $d)$ is due to the triangle inequality; $e)$ is because of \eqref{constraint}.
 Denote
$$r:=\sup_{t\geq 0}\left\{(\alpha+\theta)\cdot\|x(t)-x^*\|+\frac{\|\dot{x}(t)\|}{2}\right\}.$$ Since $\{x(t)\}_{t\geq 0}$  and  $\{\dot{x}(t)\}_{t\geq 0}$ are both bounded, we have $r<+\infty$.
Using \eqref{rela}, we have
 \begin{align}\label{system-t2}
 &\xi(t)^2=\big(f(x(t))-f(x^*)+\frac{1}{2}\|\dot{x}(t)\|^2\big)^2\nonumber\\
 &\quad\quad\leq\left((\alpha+\theta)\cdot\|x(t)-x^*\|\cdot\|\dot{x}(t)\|+\frac{1}{2}\|\dot{x}(t)\|^2\right)^2\nonumber\\
 &\quad\quad\leq\left(\left((\alpha+\theta)\cdot\|x(t)-x^*\|+\frac{1}{2}\|\dot{x}(t)\|\right)\cdot\|\dot{x}(t)\|\right)^2\nonumber\\
 &\quad\quad \leq\ r^2\|\dot{x}(t)\|^2.
 \end{align}
Combining (\ref{system-t1}) and (\ref{system-t2}), we have $\xi(t)^2\leq -\frac{r^2}{\alpha}\dot{\xi}(t)$, or equivalently,
   \begin{align}\label{moto}
 -\frac{\alpha}{r^2}d t\leq\frac{d\xi}{\xi^2}.
 \end{align}
 Taking the integral of both sides from $0$ to $t$ and noting that $\xi(0)\geq f(x(0))-\min f >0$ (we have assumed that $f(x(0))-\min f>0$), we get
 $$ -\frac{\alpha}{r^2}t\leq\frac{1}{\xi(0)}-\frac{1}{\xi(t)},$$ and thus $\xi(t)\leq \frac{1}{\frac{\alpha}{r^2}t+\frac{1}{\xi(0)}}$. Since $f(x(t))-\min f\leq\xi(t)$, we have
$$f(x(t))-\min f\leq \frac{1}{\frac{\alpha}{r^2}t+\frac{1}{\xi(0)}}$$
Thus we have derived the asymptotic sublinear convergence rate for $f(x(t))-\min f$.

\section{Convergence analysis of Heavy-ball}
In this section, we prove  convergence rates of Heavy-ball method. The core of the proof is to construct a proper Lyapunov function. The expression of $\xi(t)$ in (\ref{lyas}) suggests the Lyapunov function be of the form $f(x^k)+C\|x^{k+1}-x^k\|^2$ for some $C>0$. In fact, we have the following sufficient descent lemma.  \emph{All the technical proofs for the rest of the paper, will be provided in the supplementary materials}.
\begin{lemma}\label{le1}
Suppose $f$ is convex with $L$-Lipschitz gradient and  $\min f>-\infty$. Let $(x^k)_{k\geq 0}$ be generated by the Heavy-ball method with non-increasing $(\beta_k)_{k\geq 0}\subseteq [0,1)$. By choosing the step size
 $$\gamma_k=\frac{2(1-\beta_k)c}{L}$$
with fixed $0<c<1$, we have
\begin{align}\label{nresult}
   & \left[f(x^{k})+\frac{\beta_k}{2\gamma_k}\|x^{k}-x^{k-1}\|^2\right]\nonumber\\
   &\qquad -\left[f(x^{k+1})+\frac{\beta_{k+1}}{2\gamma_{k+1}}\|x^{k+1}-x^{k}\|^2\right]\nonumber\\
   &\qquad\geq\frac{(1-c)L}{2c}\|x^{k+1}-x^{k}\|^2.
\end{align}
\end{lemma}
According to Lemma \ref{le1}, a potentially useful Lyapunov function is $\left[f(x^{k})+\frac{\beta_k}{2\gamma_k}\|x^{k}-x^{k-1}\|^2\right]_{k\geq 0}$, as it has the descent property shown in (\ref{nresult}). However, it does not fulfill the relation in \eqref{system-t2} ~\footnote{That is, we are not able build a useful error relation for $f(x^{k})+\frac{\beta_k}{2\gamma_k}\|x^{k}-x^{k-1}\|^2$.}. Therefore, we rewrite (\ref{nresult}), so that the new right-hand-side contains something like $\|x^{k+1}-x^{k}\|^2+\|x^{k}-x^{k-1}\|^2$. It turns out that a better Lyapunov function reads
\begin{align}\label{Lyapunov}
   \xi_k:= f(x^k)+\delta_k\|x^{k}-x^{k-1}\|^2-\min f,
\end{align}
where
\begin{equation}\label{eq:delta}
    \delta_k:=\frac{\beta_k}{2\gamma_k}+\frac{1}{2}\left(\frac{1-\beta_k}{\gamma_k}-\frac{L}{2}\right).
\end{equation}
We can see $\xi_k$ is in line with the discretization of (\ref{lyas}).

Given the Lyapunov function in (\ref{Lyapunov}), we present a key technical lemma.
\begin{lemma}\label{lem-sub}
Suppose the assumptions of Lemma \ref{le1} hold.
 Let $\overline{x^k}$ denote the projection of $x^k$ onto $\emph{arg}\min{f}$, assumed to exist, and define
\begin{align}\label{eq:alpha}
 \varepsilon_k: =\frac{4c\delta_{k}^2}{(1-c)L}+\frac{4c}{(1-c)L\gamma^2_k}.
\end{align}
Then it holds that
\begin{align}\label{result-sub}
   (\xi_{k})^2&\leq\varepsilon_k\times(\xi_k-\xi_{k+1})\nonumber\\
   &\times(2\|x^{k}-\overline{x^{k}}\|^2+\|x^{k}-x^{k-1}\|^2).
\end{align}

\end{lemma}

We see that \eqref{result-sub} is the discretization  of \eqref{moto} if $\sup_{k}\{\varepsilon_k\cdot(2\|x^{k}-\overline{x^{k}}\|^2+\|x^{k}-x^{k-1}\|^2)\}<+\infty$.

\subsection{Sublinear convergence}
We present the non-ergodic $O(\frac{1}{k})$ convergence rate of the function value. This rate holds when $(\beta_k)_{k\geq 0}\in (0,1)$. We define
 \begin{align}\label{rrr}
 R:=\sup_{k\geq 0}\sup_{x^*\in \argmin f}\{\|x^k-x^*\|^2\}.
 \end{align}
 In our following settings, we can see   it actually holds that $R<+\infty$.
\begin{theorem}\label{sub-n}
Under the assumptions of Lemma \ref{le1} and assumptions that $0<\inf_k\beta_k\leq\beta_k\leq\beta_0<1$ and $f$ is coercive, We have
\begin{equation}\label{sub-n-result}
    f(x^k)-\min f\leq \frac{4R\cdot \sup_k\{\varepsilon_k\}}{k}.
\end{equation}
\end{theorem}

To our best knowledge, this is the first non-ergodic result established for Heavy-ball algorithm in convex case. The definition of $\varepsilon_k$ implies $ \sup_k\{\varepsilon_k\}=O(L)$, so it holds that
$$f(x^k)-\min f=O\left(\frac{R\cdot L}{k}\right),$$ which is on the same order of complexity as that in gradient descent.

The coercivity assumption on $f$ is crucial for Theorem \ref{sub-n}. When the function $f$ fails to be coercive, we need to assume summable $(\beta_k)_{k\geq0}$ instead.

\begin{corollary}\label{lem-b}
Suppose the assumptions of Lemma \ref{le1} hold, and $\sum_{k}\beta_k<+\infty$.\footnote{A classical example is $\beta_k=\frac{1}{k^{\theta}}$,  where $\theta>1$.}
Let $(x^k)_{k\geq 0}$ be generated by the  Heavy-ball algorithm and $f$ be coercive. Then,
$$f(x^k)-\min f\leq \frac{4R\cdot \sup_k\{\varepsilon_k\}}{k}.$$
\end{corollary}

\subsection{Linear convergence with restricted strong convexity}
We say the function $f$ satisfies a restricted strongly convex condition \citep{lai2013augmented}, if
\begin{equation}\label{condtion}
    f(x)-\min f\geq \nu\|x-\overline{x}\|^2,
\end{equation}
where $\overline{x}$ is the projection of $x$ onto the set $\textrm{arg}\min f$, and $\nu>0$. Restricted strong convexity is weaker than the strong convexity. For example, let us consider the function $\frac{1}{2}\|Ax-b\|^2$ with $b\in \textrm{range}(A)$. When $A$ fails to be full row-rank, $\frac{1}{2}\|Ax-b\|^2$ is not strongly convex but restricted strongly convex.
\begin{theorem}\label{thm-linear}
Suppose the assumptions of Theorem \ref{sub-n} hold, and that $f$ satisfies condition (\ref{condtion}).
 Then we have
\begin{equation*}\label{thm-linear-result}
   f(x^k)-\min f\leq \omega^k,
\end{equation*}
for $\omega^k:=\frac{\ell}{1+\ell}\in(0,1)$ and $\ell:=\sup_{k}\left\{\varepsilon_k(\frac{1}{\delta_k}+\frac{2}{\nu})\right\}$.
\end{theorem}

Our result improves the linear convergence established by \cite{ghadimi2015global} in two aspects: Firstly, The strongly convex assumption is weakened to (\ref{condtion}). Secondly, The step size and inertial parameter are chosen independent of the strongly convex constants.



\section{Cyclic coordinate descent Heavy-ball algorithm}
In this section, we consider the multi-block version of Heavy-ball algorithm and prove its convergence rates under convexity assumption. The minimization problem reads
\begin{equation}\label{multi}
    \min_{x_1,x_2,\ldots,x_m}f(x_1,x_2,\ldots,x_m).
\end{equation}
The function $f$ is assumed to satisfy
\begin{align}\label{picl}
   \|\nabla_i f(x)-\nabla_i f(y)\|\leq L_i\|x-y\|.
\end{align}
 With (\ref{picl}), we can easily obtain
 \begin{align}\label{lpm}
   & f(x_1,x_2,\ldots,x_i^1,\ldots,x_m)\leq  f(x_1,x_2,\ldots,x_i^2,\ldots,x_m)\nonumber\\
    &\quad\quad+\langle \nabla_i f(x_1,x_2,\ldots,x_i^2,\ldots,x_m),x_i^1-x_i^2\rangle\nonumber\\
    &\quad\quad +\frac{L_i}{2}\|x_i^1-x_i^2\|^2.
 \end{align}
The proof is similar to [Lemma 1.2.3,\cite{nesterov2013introductory}], and we shall skip it here. We denote
\begin{align*}
\nabla_i^k f:=\nabla_i f(x_1^{k+1},\ldots,x_{i-1}^{k+1},x_{i}^{k}\ldots,x_m^k),\\ \quad x^k:=(x^k_1,x^k_2,\ldots,x^k_m),~\quad L:=\sum_{i=1}^m L_i,
\end{align*}
with the convention $x^{k+1}_0=x^k_1$.
 The cyclic coordinate descent inertial algorithm iterates: for $i$ from $1$ to $m$,
\begin{equation}\label{schemec}
    x^{k+1}_i=x^k_i-\gamma_{k,i}\nabla_i^k f+\beta_{k,i}(x^k_i-x^{k-1}_{i}),
\end{equation}
where $\gamma_{k,i},\beta_{k,i}>0$.
Our analysis relies on the following assumption:

\textbf{A1}: for any $i\in[1,2,\ldots,m]$, the parameters $(\beta_{k,i})_{k\geq 0}\subseteq [0,1)$ is non-increasing.

\begin{lemma}\label{md-le1}
 Let $f$ be a convex function  satisfying (\ref{picl}), and finite $\min f$. Let $(x^k)_{k\geq 0}$ be generated by scheme (\ref{schemec}) and Assumption \textbf{A1} hold. Choosing the step size
 $$\gamma_{k,i}=\frac{2(1-\beta_{k,i})c}{L_i},~\quad i\in[1,2,\ldots,m]$$
 for arbitrary fixed $0<c<1$, we have
\begin{align}\label{md-result}
    &\left[f(x^{k})+\sum_{i=1}^m\frac{\beta_{k,i}}{2\gamma_{k,i}}\|x^{k}_i-x^{k-1}_i\|^2\right]\nonumber\\
    &\quad\quad-\left[f(x^{k+1})+\sum_{i=1}^m\frac{\beta_{k+1,i}}{2\gamma_{k+1,i}}\|x^{k+1}_i-x^{k}_i\|^2\right]\nonumber\\
    &\quad\quad\geq\frac{(1-c)\underline{L}}{2c}\|x^{k+1}-x^k\|^2.
\end{align}
where $\underline{L}=\min_{i\in[1,2,\ldots,m]}\{L_i\}$.
\end{lemma}

We consider the following similar Lyapunov function in the analysis of cyclic coordinate descent Heavy-ball algorithm
\begin{align}\label{Lyapunov-m}
   \hat{\xi}_k:= f(x^k)+\sum_{i=1}^m \delta_{k,i}\|x^{k}_i-x^{k-1}_i\|^2-\min f,
\end{align}
where
\begin{equation}\label{eq:delta-m}
    \delta_{k,i}:=\frac{\beta_{k,i}}{2\gamma_{k,i}}+\frac{1}{2}\left(\frac{1-\beta_{k,i}}{\gamma_{k,i}}-\frac{L_i}{2}\right).
\end{equation}
Then we have the following lemma.
\begin{lemma}\label{lem-sub-m}
Suppose the conditions Lemma \ref{md-le1} hold.
 Let $\overline{x^k}$ denote the projection of $x^k$ onto $\emph{arg}\min{f}$, assumed to exist, and define
\begin{align}
 \hat{\varepsilon}_k&:=\max\left\{\frac{4c\cdot\sum_{i=1}^m\left(\delta_{k+1,i}^2+\frac{1}{\gamma^2_{k,i}}\right)}{(1-c)\underline{L}},\frac{4c\cdot m\cdot L}{(1-c)\underline{L}}\right\}.
\end{align}
It holds  that
\begin{align}\label{result-sub-m}
   (\hat{\xi}_{k})^2&\leq\hat{\varepsilon}_k(\hat{\xi}_k-\hat{\xi}_{k+1})\nonumber\\
   &\times(2\|x^{k}-\overline{x^{k}}\|^2+\|x^{k}-x^{k-1}\|^2).
\end{align}
\end{lemma}

\subsection{Sublinear convergence of  cyclic coordinate descent Heavy-ball algorithm}
We show the $O(1/k)$ convergence rate of cyclic coordinate descent Heavy-ball algorithm for coercive $f$.

\begin{theorem}\label{them-sub-m}
Suppose the conditions of Lemma \ref{md-le1} hold, $f$ is coercive, and
 $$0<\inf_k\beta_{k,i}\leq\beta_{k,i}\leq\beta_0<1,~i\in[1,2,\ldots,m].$$
Then we have
\begin{align}
   f(x^k)-\min f= O\left(\frac{4R\cdot \sup_{k}\{\hat{\varepsilon}_k\}}{k}\right),
\end{align}
where $R$ is given by \eqref{rrr}.
\end{theorem}
We readily check that $\sup_{k}\{\hat{\varepsilon}_k\}=O(mL)$, where $m$ is number of the block. Therefore, the cyclic inertial algorithm converges with the rate $O\left(\frac{m\cdot R\cdot L}{k}\right)$. Compared with the results in \citep{sun2015improved}, this rate is on the same order as that of cyclic block coordinate descent in general convex setting.

\subsection{Linear convergence of cyclic coordinate descent Heavy-ball algorithm}

Under the same assumption of restricted strong convexity, we derive the linear convergence rate for cyclic coordinate descent Heavy-ball algorithm.
\begin{theorem}\label{them-linear-m}
Suppose the conditions of Lemma \ref{md-le1} hold, $f$ satisfies (\ref{condtion}), and
 $$0<\inf_k\beta_{k,i}\leq\beta_{k,i}\leq\beta_0<1,~i\in[1,2,\ldots,m].$$
Then we have
\begin{align}
   f(x^k)-\min f\leq (\hat{\omega})^k
\end{align}
for some $\hat{\omega}=\frac{\hat{\ell}}{1+\hat{\ell}}\in(0,1)$, and $\hat{\ell}:=\sup_{k}\big\{\hat{\varepsilon}_k+\frac{2}{\nu}+\frac{1}{\min_i\{\delta_{k,i}\}}\big\}$.
\end{theorem}

This result can be extended to the essentially cyclic Heavy-ball algorithm. The essentially cyclic index selection strategy \citep{sun2017asynchronous}, which generalizes the cyclic update rule, is defined as follows: there is an $M\in\mathbb{N}$, $M\ge m$, such that each block $i\in\{1,2,\ldots,m\}$ is updated at least once in a window of $M$.

\section{Stochastic coordinate descent Heavy-ball algorithm}
For the stochastic index selection strategy, in the $k$-th iteration, we pick $i_k$ uniformly from $[1,2,\ldots,m]$ and iterate
\begin{align}\label{schemes}
   \left\{\begin{array}{c}
            x^{k+1}_{i_k}=x^k_{i_k}-\gamma_{k}\nabla_{i_k} f(x^k)+\beta_{k}(x^k_{i_k}-x^{k-1}_{i_k}), \\
            x^{k+1}_i=x^k_i,~\textrm{if}~i\neq i_k.
          \end{array}
   \right.
\end{align}
In this section, we make the following assumption

\textbf{A2}:  the parameters $(\beta_{k})_{k\geq 0}\subseteq [0,\sqrt{m})$ is non-increasing.

Assumption  \textbf{A2} is quite different from previous requirement that $(\beta_{k})_{k\geq 0}$ which are constrained on $[0,1)$. This difference comes from the uniformly stochastic selection of the index.

\begin{lemma}\label{ms-le1}
 Let $f$ be a convex function  whose gradient is Lipschitz continuous with $L$, and finite $\min f$. Let $(x^k)_{k\geq 0}$ be generated by scheme (\ref{schemes}) and Assumption \textbf{A2} be satisfied. Choose the step size
 $$\gamma_{k}=\frac{2(1-\beta_{k}/\sqrt{m})c}{L}$$
 for arbitrary fixed $0<c<1$. Then, we can obtain
\begin{align}\label{ms-result}
&\left[\EE f(x^{k})+\frac{\beta_k}{2\sqrt{m}\gamma_k}\EE\|x^k-x^{k-1}\|^2\right]\nonumber\\
&\quad\quad-\left[\EE f(x^{k+1})+\frac{\beta_{k+1}}{2\sqrt{m}\gamma_{k+1}}\EE\|x^{k+1}-x^k\|^2\right]\nonumber\\
&\quad\quad\geq\frac{(1-c)L}{2c}\EE\|x^{k+1}-x^k\|^2.
\end{align}
\end{lemma}

Similarly, we consider the following function
\begin{align}\label{Lyapunov-s}
   \bar{\xi}_k:= f(x^k)+ \bar{\delta}_{k}\|x^k-x^{k-1}\|^2-\min f,
\end{align}
where
\begin{equation}\label{eq:delta-s}
  \bar{\delta}_{k}:=\frac{\beta_{k}}{2\sqrt{m}\gamma_{k}}+\frac{1}{2}\left(\frac{1-\beta_{k}/\sqrt{m}}{\gamma_{k}}-\frac{L}{2}\right).
\end{equation}
Different from the previous analyses, the Lyapunov function considered here is $\EE \bar{\xi}_k$ instead of $\bar{\xi}_k$. Naturally, the sufficient descent property is established in the sense of expectation.

\begin{lemma}\label{lem-sub-s}
Suppose the conditions of Lemma \ref{ms-le1} hold.
 Let $\overline{x^k}$ denote the projection of $x^k$ onto $\emph{arg}\min{f}$, assumed to exist, and define
\begin{align}\label{eq:alpha-s}
 \bar{\varepsilon}_k: =\frac{4c\delta_{k}^2}{(1-c)L}+\frac{8cm}{(1-c)L\gamma^2_k}.
\end{align}
Then it holds
\begin{align}\label{result-sub2}
   (\EE\bar{\xi}_{k})^2&\leq\bar{\varepsilon}_k\cdot(\EE\bar{\xi}_k-\EE\bar{\xi}_{k+1})\nonumber\\
   &\times(\EE\|x^{k}-\overline{x^{k}}\|^2+\EE\|x^{k}-x^{k-1}\|^2).
\end{align}
\end{lemma}

\subsection{Sublinear convergence of stochastic coordinate descent Heavy-ball algorithm}
Due to that the sufficient descent condition involves expectations, even using the coercivity of $f$, we cannot obtain the boundedness of the generated points. Therefore, we  first present a result by assuming the smoothness of $f$ only.
\begin{theorem}\label{ms-th1}
Suppose that the assumptions of Lemma \ref{ms-le1} hold.
Then we have
\begin{equation}\label{ms-th-result}
    \min_{0\leq i\leq k}\EE\|\nabla f(x^k)\|=o\left(\frac{1}{\sqrt{k}}\right).
\end{equation}
\end{theorem}
We remark that Theorem \eqref{ms-th1} also holds for nonconvex functions.
To obtain the sublinear convergence rate on the function values, we need a boundedness assumption. Precisely, the assumption is

\textbf{A3}:  the sequence $(x^k)_{k\geq 0}$ satisfies
$$\bar{R}:=\sup_{k}\left\{\EE\|x^k-\overline{x^k}\|^2\right\}<+\infty.$$

Under assumption \textbf{A3}, we are able to show the non-ergodic convergence sublinear convergence rates of the expected objective values.
\begin{theorem}\label{ms-th1-f}
Suppose that the assumptions of Lemma \ref{ms-le1} and \textbf{A3} hold.
Then we have
\begin{equation}
    \EE f(x^k)-\min f=O\left(\frac{4\bar{R}\cdot\sup_{k}\{\bar{\varepsilon}_k\}}{k}\right).
\end{equation}
\end{theorem}
\begin{figure*}[!htb]
    \centering
  \subfigure[]{ \includegraphics[width=0.24\textwidth]{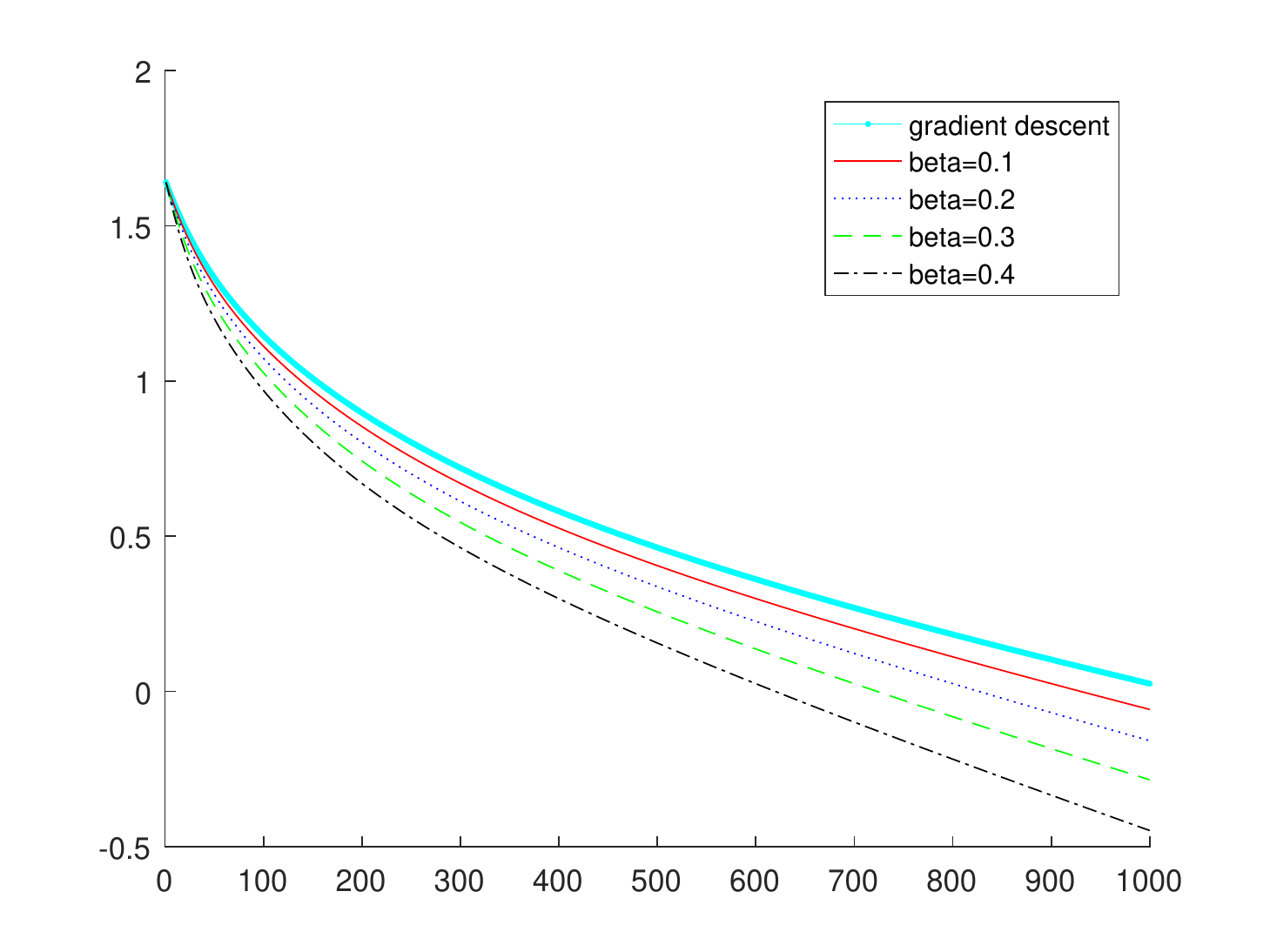}}    \subfigure[]{ \includegraphics[width=0.24\textwidth]{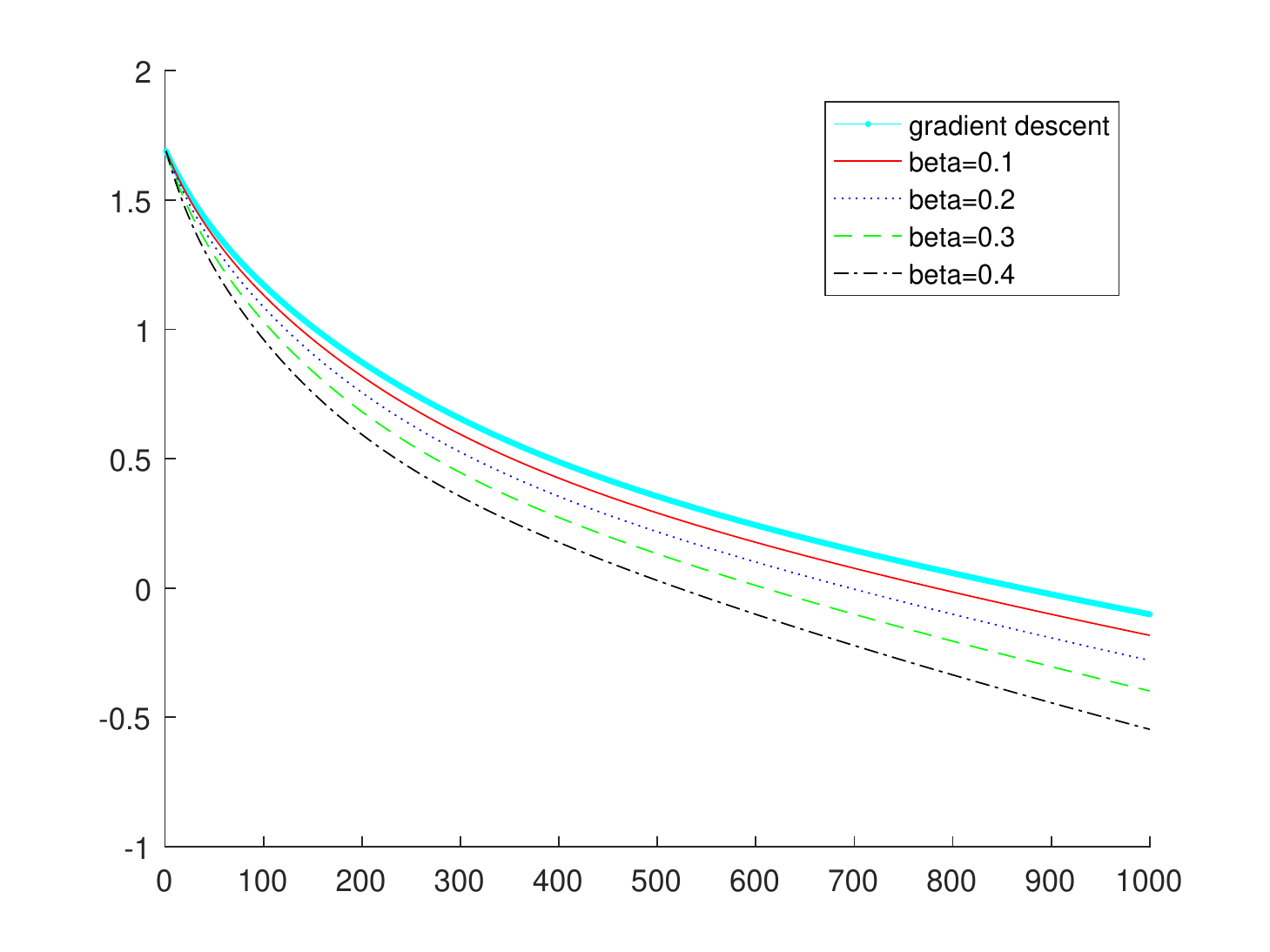}}
    \subfigure[]{ \includegraphics[width=0.24\textwidth]{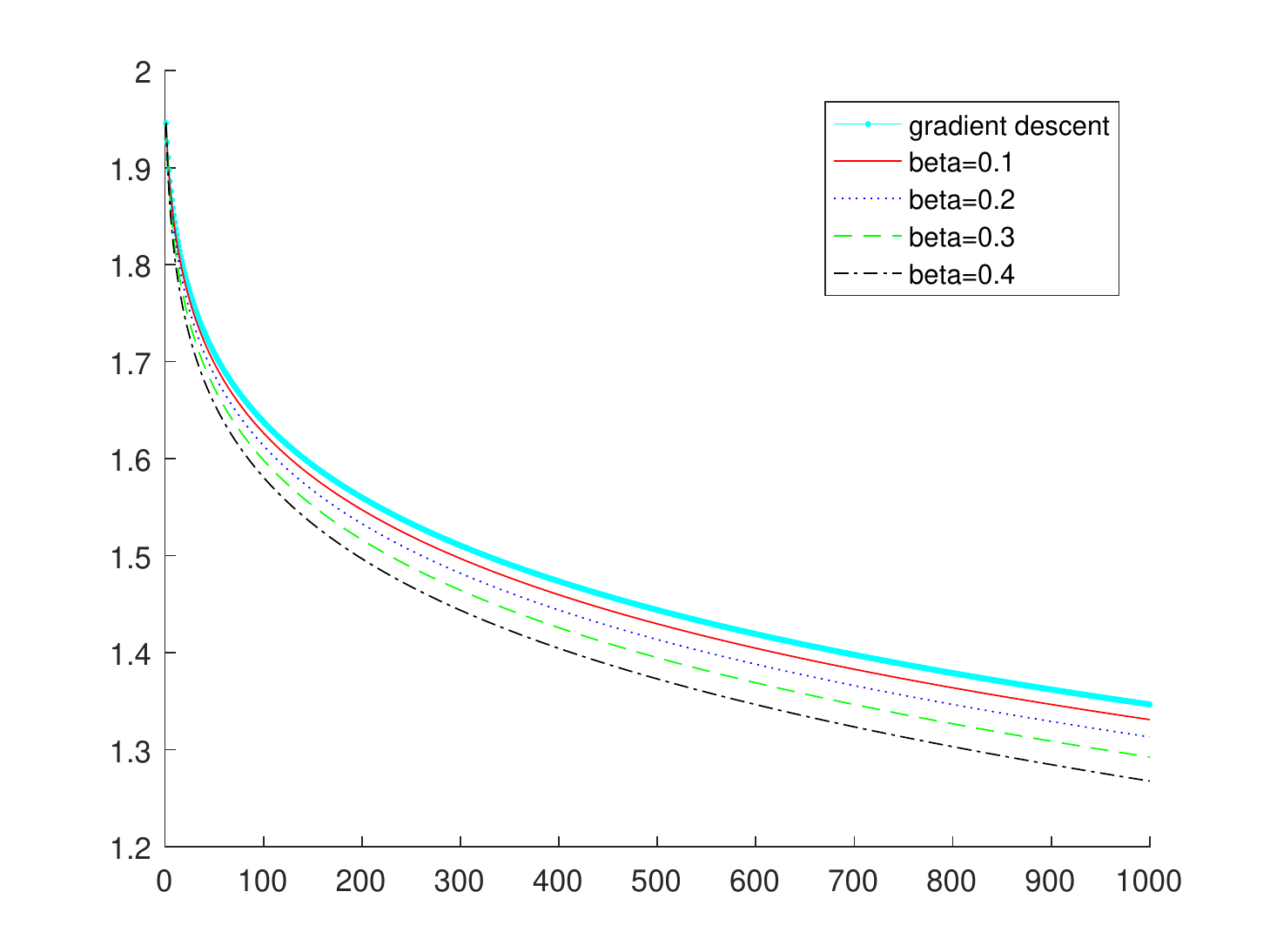}}    \subfigure[]{ \includegraphics[width=0.24\textwidth]{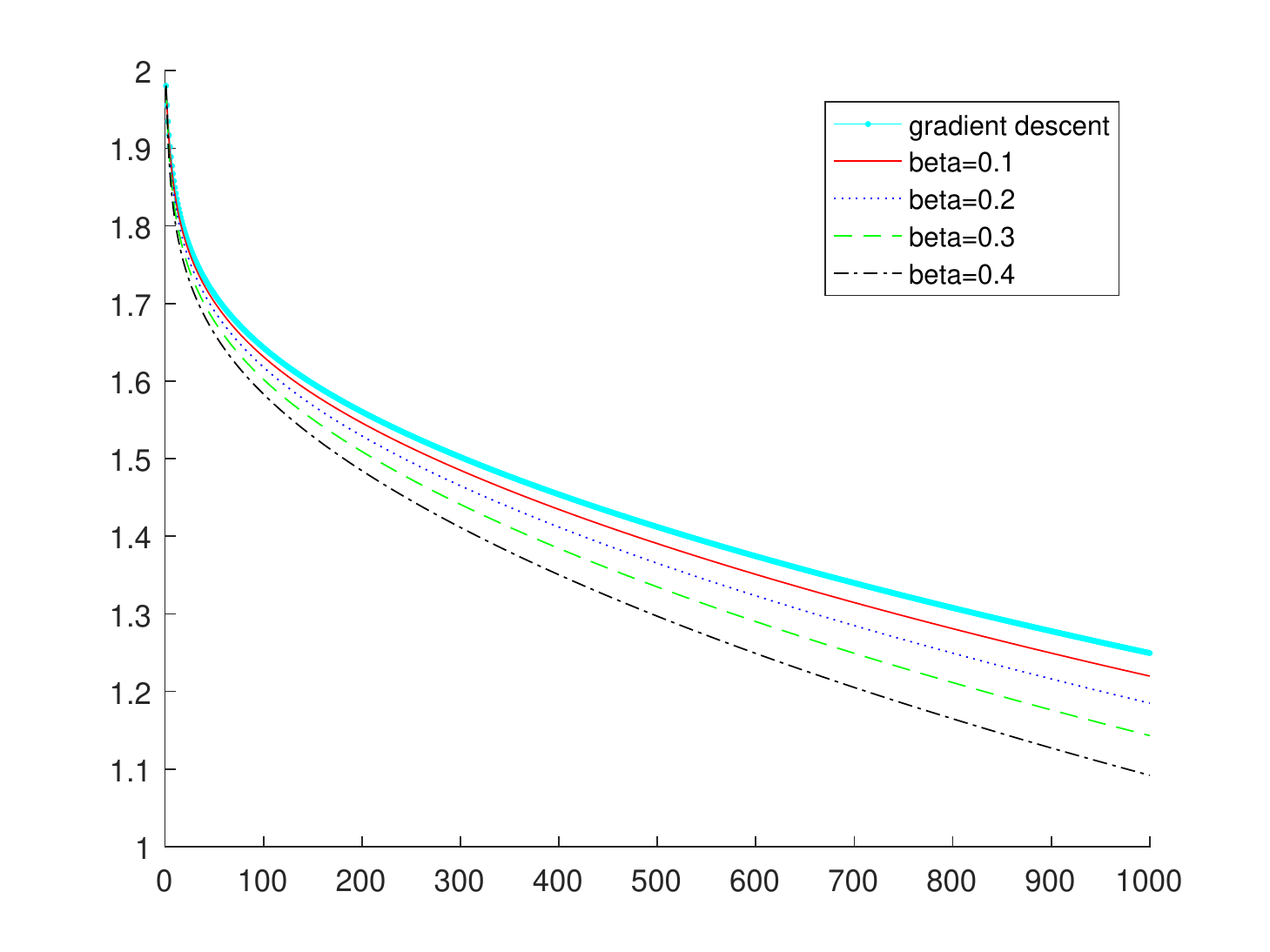}}\\
  \subfigure[]{\includegraphics[width=0.24\textwidth]{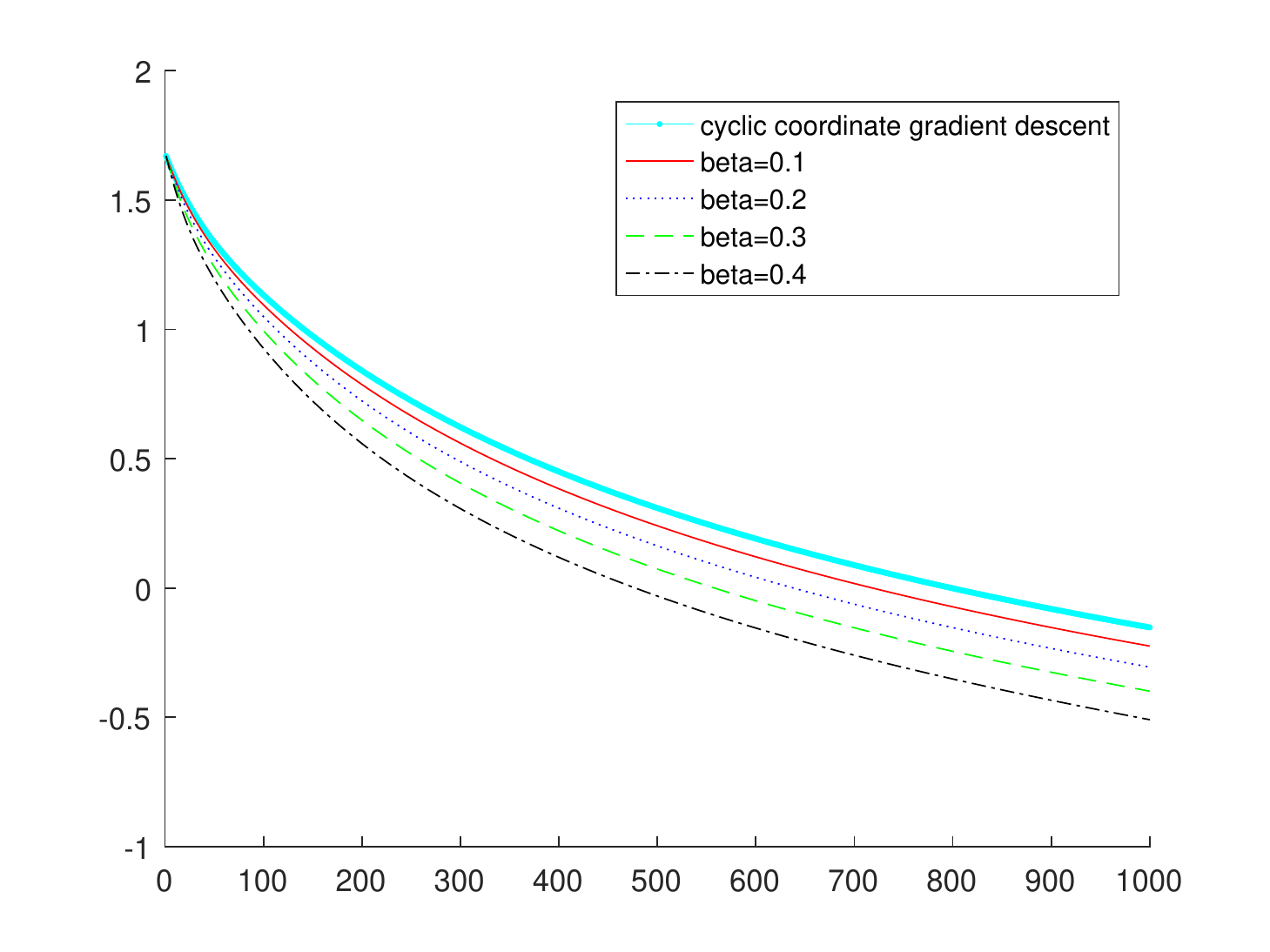}}
     \subfigure[]{ \includegraphics[width=0.24\textwidth]{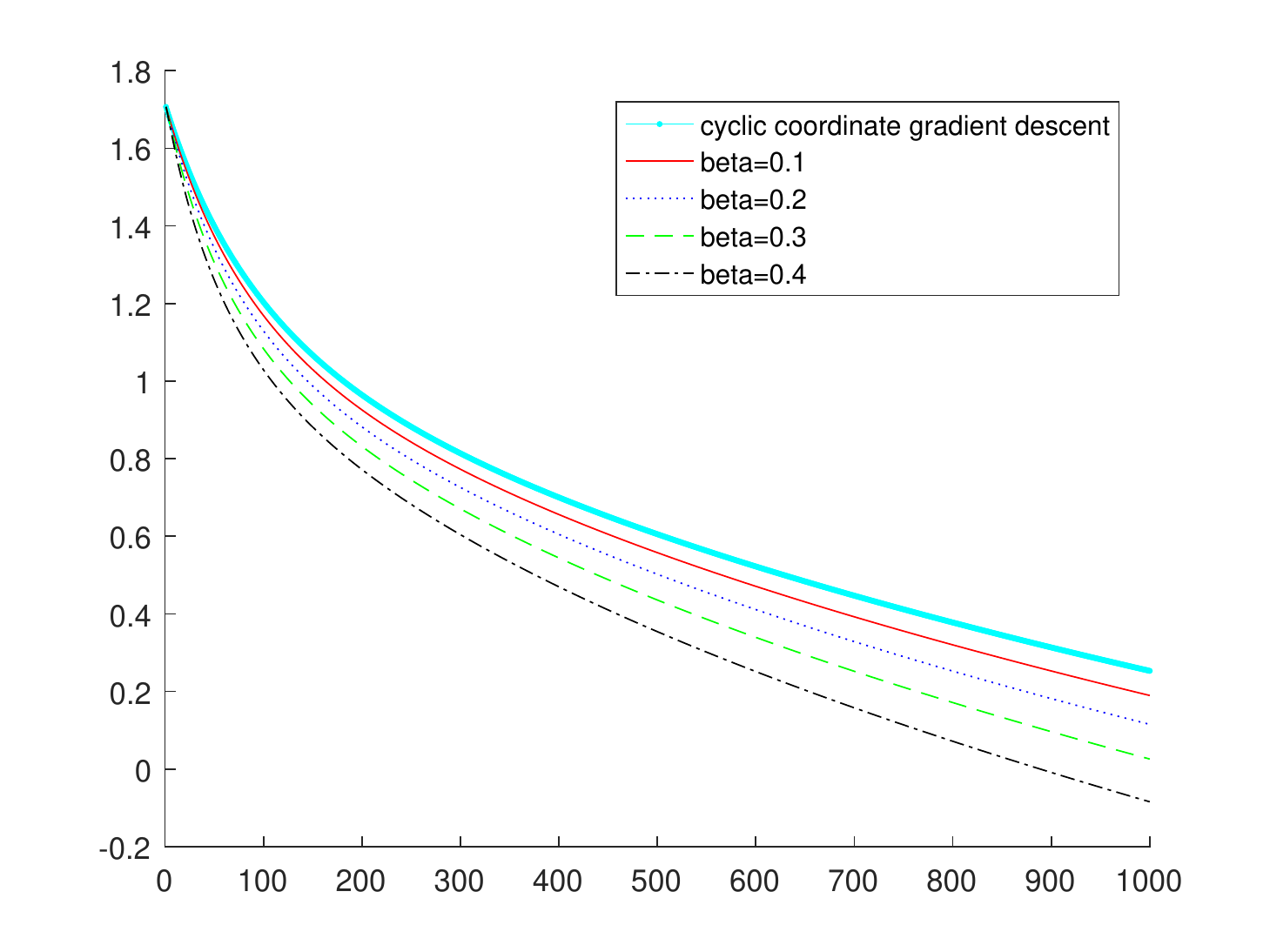}}
      \subfigure[]{\includegraphics[width=0.24\textwidth]{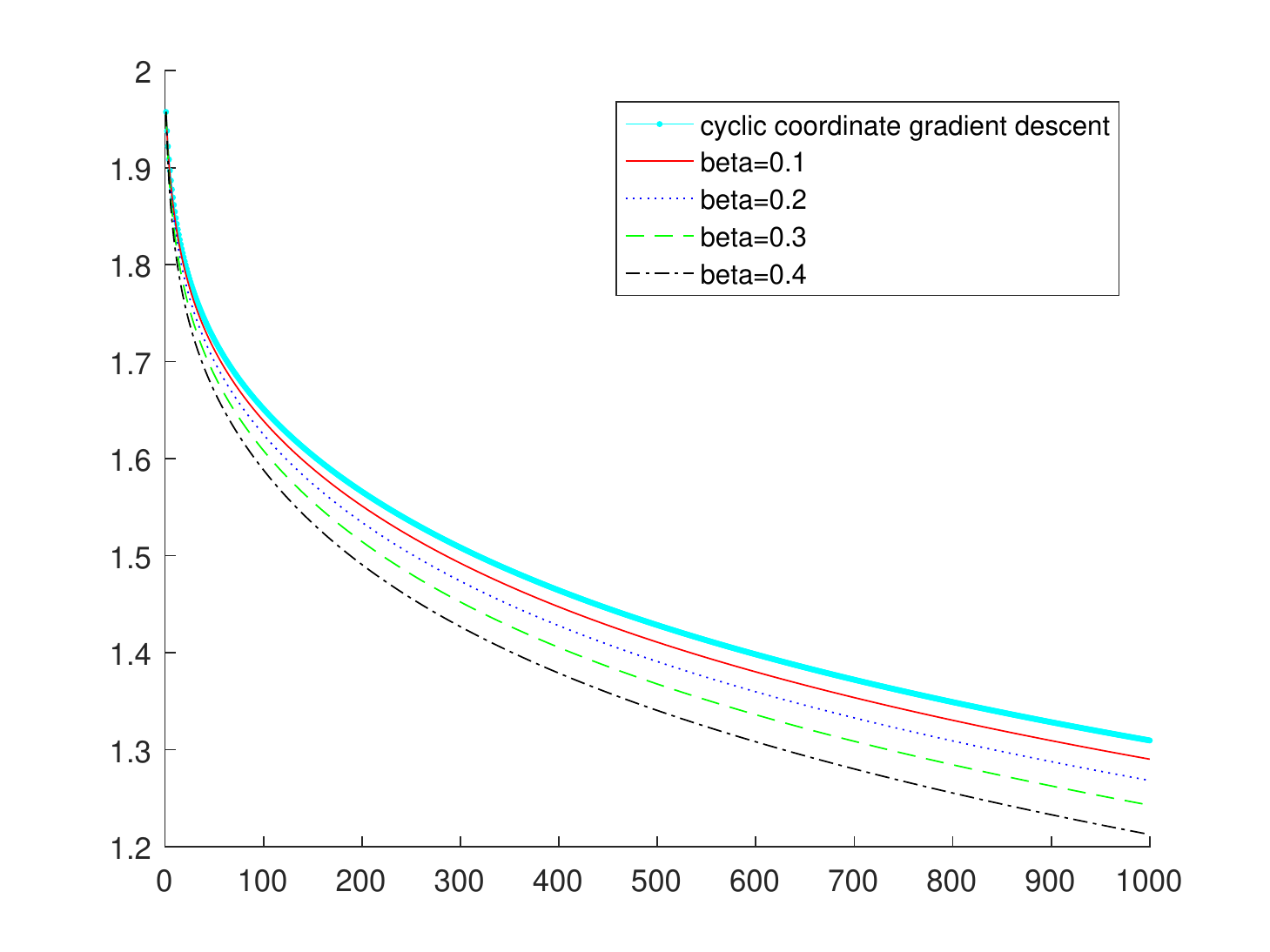}}
     \subfigure[]{ \includegraphics[width=0.24\textwidth]{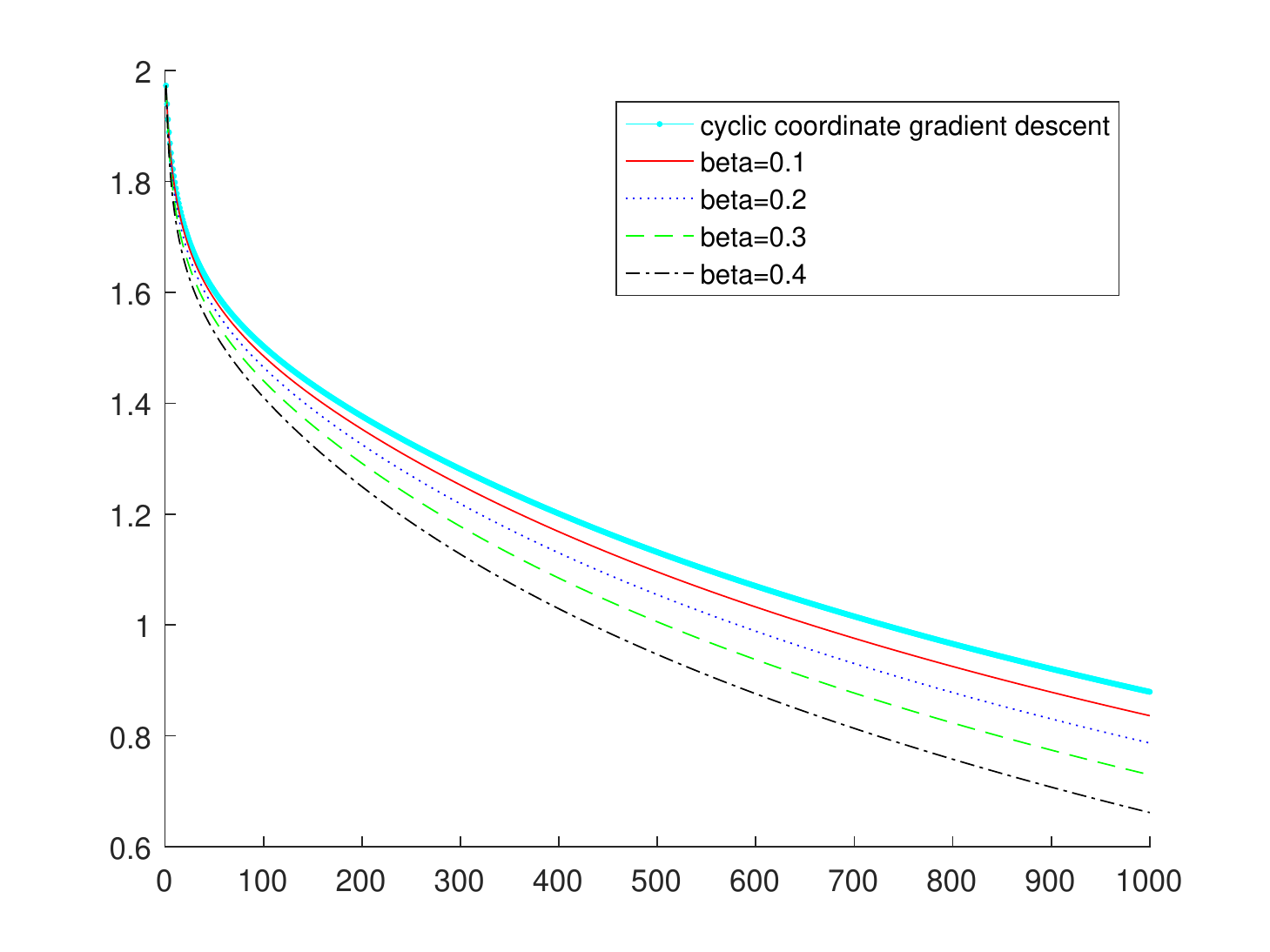}}\\
  \subfigure[]{ \includegraphics[width=0.24\textwidth]{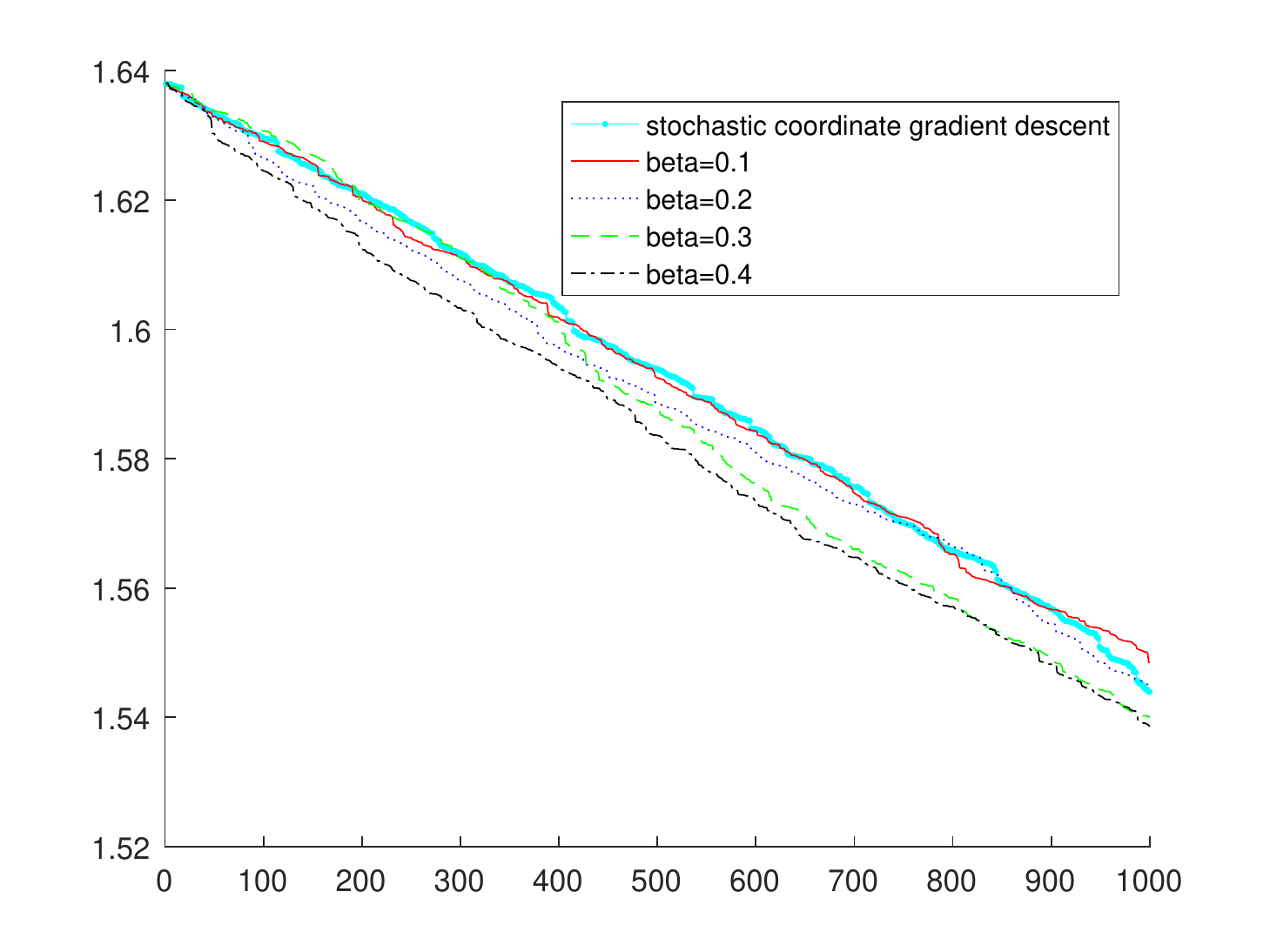}}
     \subfigure[]{ \includegraphics[width=0.24\textwidth]{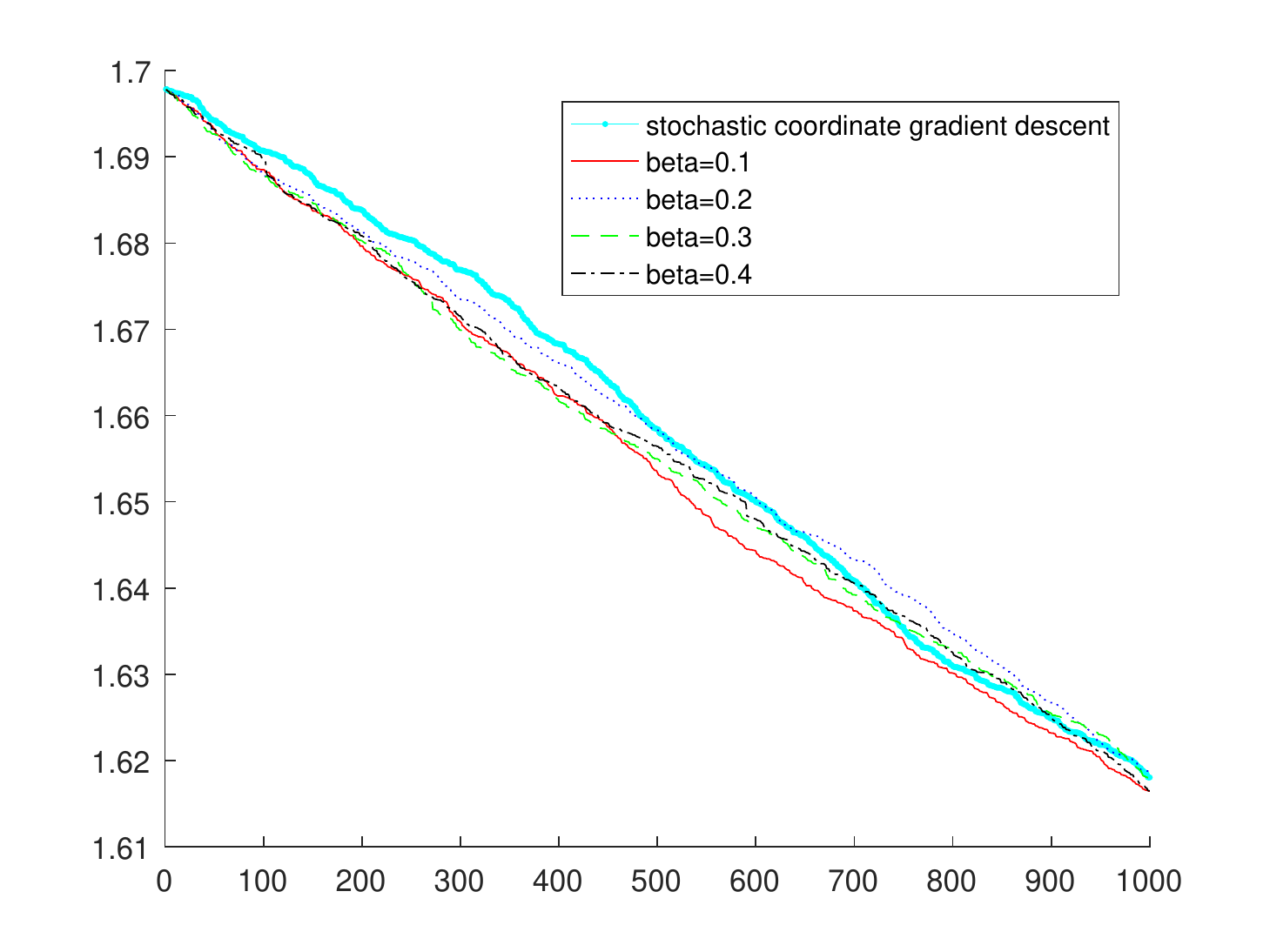}}
       \subfigure[]{ \includegraphics[width=0.24\textwidth]{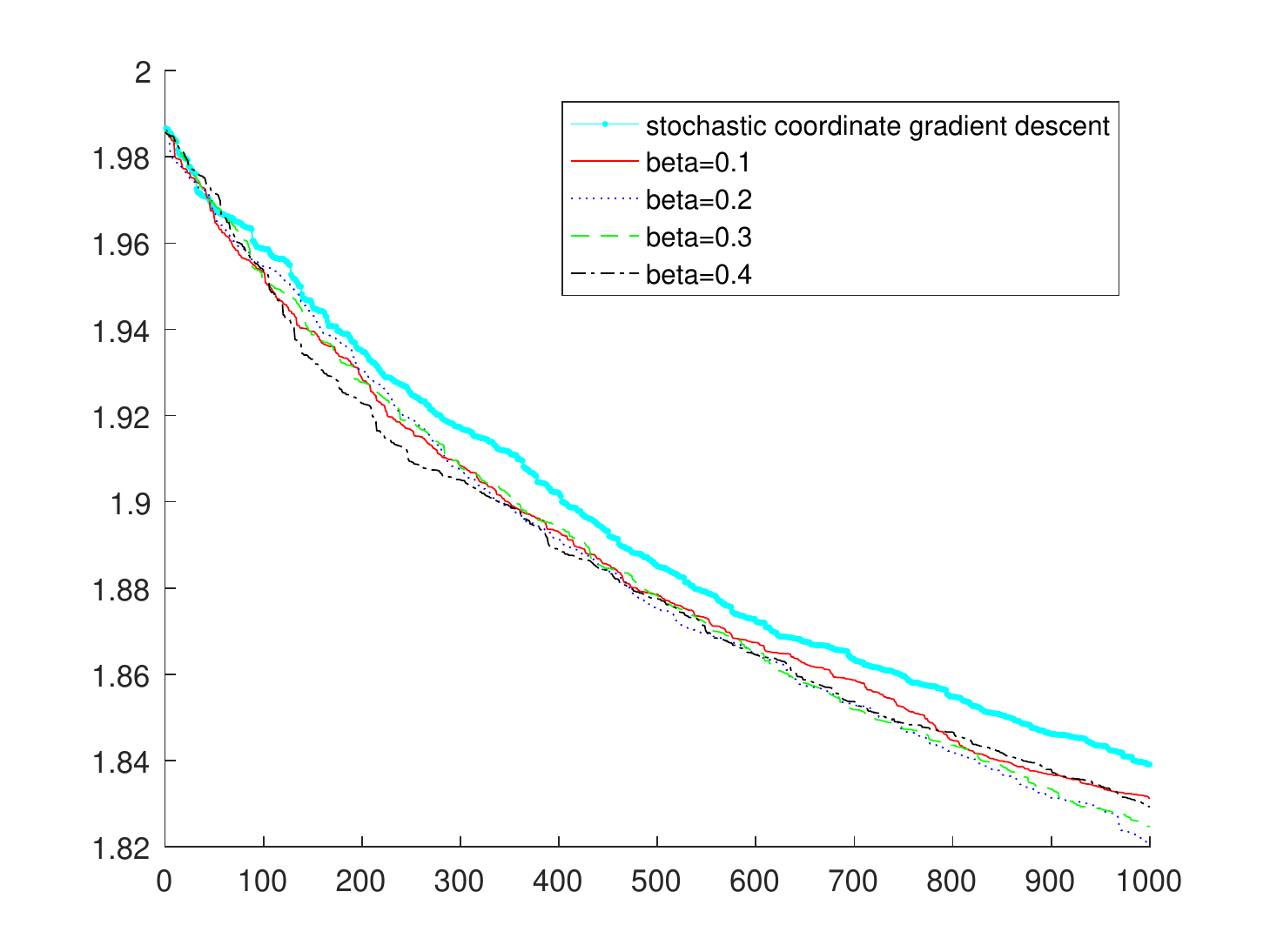}}
     \subfigure[]{ \includegraphics[width=0.24\textwidth]{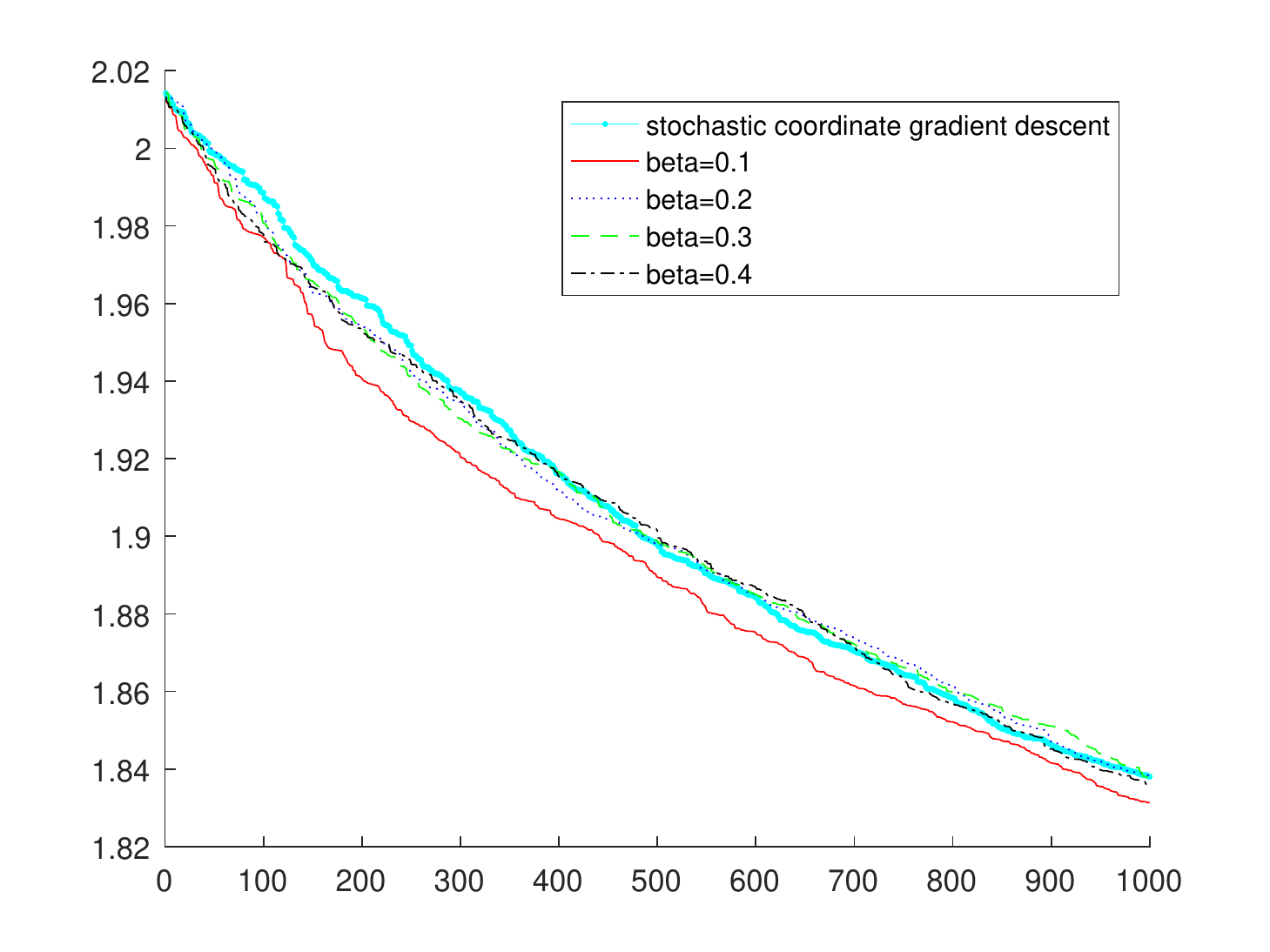}}
    \caption{Objective values of inertial algorithms v.s. number of epochs for linear regression and logistic regression tasks. Panels (a)-(d): Heavy-ball algorithms comparisons  for linear regression with Gaussian data (a); for  linear regression with Bernoulli data (b);   for   logistic regression with Gaussian data (c); for   logistic regression with Bernoulli data (d).
Panels (e)-(h): Cyclic coordinate descent  Heavy-ball algorithms comparisons  for linear regression with Gaussian data (e); for  linear regression with Bernoulli data (f); for logistic regression with Gaussian data (g); for logistic regression with Bernoulli data (h).
Panels (i)-(l): Stochastic coordinate descent  Heavy-ball algorithms comparisons  for linear regression with Gaussian data (i); for  linear regression with Bernoulli data (j);   for   logistic regression with Gaussian data (k); for logistic regression with Bernoulli data (l).}
\end{figure*}

\subsection{Linear convergence of  stochastic coordinate descent Heavy-ball algorithm}
The linear convergence rate of stochastic coordinate descent Heavy-ball algorithm is similar to previous ones. By assuming the restricted strongly convex condition, the linear convergence rate of the expected objective values can be proved.
\begin{theorem}\label{th:linear-convergence}
Suppose that the assumptions in Lemma \ref{ms-le1} hold, and the function satisfies the restricted strongly convex condition (\ref{condtion}). Let  $(x^k)_{k\geq 0}$ be generated by the scheme (\ref{schemes}). Then we have
\begin{align}
\EE f(x^k)-\min f\leq (\bar{w})^k,
\end{align}
where $\bar{w}:=\frac{\bar{\ell}}{1+\bar{\ell}}\in(0,1)$, and $\bar{\ell}:=\sup_{k}\{\bar{\varepsilon}_k+\frac{1}{\nu}+\frac{1}{\bar{\delta}_{k}}\}$.
\end{theorem}

While we only consider the uniform probability selection strategy here, the same convergence results can be easily extended to the non-uniform  probability selection strategy.

\section{Applications to decentralized optimization}
We apply the analysis to the following decentralized optimization problem
\begin{align*}
\min_{x\in\textrm{R}^{ n}}\left\{\sum_{i=1}^m f_i(x)\right\},
\end{align*}
where $f_i$ is differentiable and $\nabla f_i$ is $L_i$-Lipschitz.
Denote by $x(i)\in \mathbb{R}^n$ the local copy of $x$ at node $i$  and $X:=(x(1),x(2),\ldots,x(m))^{\top}$.
In the community of decentralized algorithms, rather than directly solving the problem, following penalty formulation instead has been proposed
\begin{align}\label{dectral}
\min_{X\in \mathbb{R}^{m\times n}}\left\{F(X)=f(X)+\frac{X^{\top}(I-W)X}{2\alpha}\right\},
\end{align}
where $W=(w_{i,j})\in \mathbb{R}^{m\times m}$ is the mixing matrix, and $f(X):=\sum_{i=1}^m f_i(x(i))$, and $I$ is the unit matrix. It is easy to see that $\nabla F$ is Lipschitz with the constant $L_F:=\max_{i}\{L_i\}+\frac{1-\lambda_{\min}(W)}{\alpha}$, here $\lambda_{\min}(W)$ is minimum eigenvalue of $W$.
Researchers consider the decentralized gradient descent (DGD) \citep{nedic2009distributed}, which is essentially the gradient descent applied to \eqref{dectral} with stepsize being equal to $\alpha$. This algorithm can be implemented over a connected network, in which the agents communicate with their neighbors and make full use of the computing resources of all nodes.
 Alternatively, we can use the Heavy-ball method by choosing the stepsize $\alpha$, that is,
 \begin{align*}
X^{k+1}&=X^{k}-\alpha \nabla F(X^k)+\beta(X^{k}-X^{k-1})\\
&=WX^k-\alpha \nabla f(X^k)+\beta(X^{k}-X^{k-1}).
\end{align*}
 For node $i$,   the local scheme is then
\begin{align*}
x^{k+1}(i)&=\sum_{j\in \mathcal{N}(i)} w_{i,j}  x^{k}(j)-\alpha\nabla f_i(x^k(i))\\
&+\beta (x^k(i)-x^{k-1}(i)),
\end{align*}
where $x^k(i)$ is the copy of the variable $x^k$ in node $i$ in the $k$th iteration and $\mathcal{N}(i)$ denotes the neighbors of node $i$. In the global scheme, it is basically Heavy-ball algorithm. Thus, we can apply our theoretical findings to this algorithm. To guarantee the convergence, we just need
$$0<\alpha<\frac{2(1-\beta)}{L_F}=\frac{2(1-\beta)}{\max_{i}\{L_i\}+\frac{1-\lambda_{\min}(W)}{\alpha}}.$$
After simplification, we then get
$$\alpha\cdot\max_{i}\{L_i\}<1+\lambda_{\min}(W)-2\beta.$$
In a word, we need the requirements $$0\leq \beta<\frac{1+\lambda_{\min}(W)}{2},\,~0<\alpha<\frac{1-2\beta+\lambda_{\min}(W)}{\max_{i}\{L_i\}}.$$ The convergence result for decentralized Heavy-ball method directly follows from our previous theoretical findings and  can be summarized as below.
\begin{corollary}
Assume that $f_i$ is convex and differentiable,  and   $\nabla f_i$ is Lipschitz with $L_i$. Let $0\leq \beta<\frac{1+\lambda_{\min}(W)}{2}$, and the sequence $(X^k)_{k\geq 0}$ be generated
by the decentralized Heavy-ball method. For any fixed stepsize $0<\alpha<\frac{1-2\beta+\lambda_{\min}(W)}{\max_{i}\{L_i\}}$, we have
\begin{align}
F(X^k)-\min F=O\left(\frac{1}{k}\right).
\end{align}
\end{corollary}

This justifies the superiority of our non-ergodic analysis.
As aforementioned in the introduction, all the existing convergence results are about the sequence $\left\{F\left(\frac{\sum_{i=1}^k X^i}{k}\right)-\min F\right\}$. However, for decentralized Heavy-ball algorithm, it is meaningless to discuss the ergodic rates, because the nodes only communicate with their neighbors.
However, our results, in this case, still hold.

\section{Experimental results}
We report the numerical simulations of Heavy-ball method applied to the linear regression problem
\begin{align}\label{LSP}
\min_{x\in\mathbb{R}^{n}}
\left\{\frac{1}{2}\sum_{i=1}^m(y_i-A_i^{\top}x)^2\right\},
\end{align}
and  the logistic regression problem
\begin{align}\label{LRP}
\min_{x\in\mathbb{R}^{n}}\left\{ \sum_{i=1}^m \log\left(1+\textrm{exp}(-y_i A_i^{\top}x)\right)+\frac{\lambda}{2}\|x\|^2\right\},
\end{align}
where $(A_i,y_i)\in \mathbb{R}^{n}\oplus \mathbb{R}$, $i=1,2,\ldots,m$.
All experiments were performed using MATLAB on an desktop with an Intel 3.4 GHz CPU. We tested the three Heavy-ball algorithms with different inertial parameters. We fixed the stepsize  as $\gamma=\frac{1}{L}$ in all numerical tests. For the stepsize, we need $2(1-\beta_k)>1$, i.e., $0\leq\beta_k<0.5$. Therefore, inertial parameters are set to
$\beta_k\equiv\beta=0,0.1,0.2,0.3,0.4$. For linear regression problem, $L=\lambda_{\max}(\sum_{i=1}^mA_i^{\top}A_i)$, where $\lambda_{\max}(\cdot)$ denotes the largest eigenvalue of a matrix; whereas for logistic regression, we have $L=\lambda_{\max}(\sum_{i=1}^mA_i^{\top}A_i)+\lambda$.    With schemes of the algorithms, for cyclic coordinate gradient descent, the function values are recorded after the whole epoch is updated; while for stochastic  coordinate gradient descent, functions values are updated after per iteration.  The special case $\beta=0$ corresponds to the gradient descent, or cyclic coordinate gradient descent, or stochastic coordinate gradient descent.  And we set $n=100$ and $m=150$. The data $A_i$ and $y_i$ were generated by the Gaussian random and Bernoulli random distributions, respectively. The maximum number of iterations was set to $1000$. For logistic regression, we set $\lambda=10^{-3}$.  We tested the three Heavy-ball algorithms for both two regression tasks with Gaussian and Bernoulli data.

As illustrated by Figure 1, larger $\beta$ leads to faster convergence for  both Heavy-ball algorithm and cyclic coordinate descent algorithm  when $\beta\in[0,0.4]$. However, for the stochastic block coordinate descent scheme, the inertial method  helps insignifically. This is  because for the stochastic case, in the $k$th iteration, the inertial terms contribute  only when $i_k=i_{k-1}$. This case, however, happens with probability $N\cdot\frac{1}{N^2}=\frac{1}{N}$ and $N$ is the number of the blocks; as $N$ is large, $i_k=i_{k-1}$ happens at low probability for just one iteration, let alone the whole iterations.   Therefore, the inertial method is actually  inactive at most iterations for stochastic block coordinate descent scheme.

To improve the practical performance of stochastic block coordinate descent Heavy-ball algorithm, another inertial scheme proposed in \cite{xu2013block} can be  recruited, in which,   a new storage $y^k$ is used. In each iteration, the algorithm employs $\gamma_k(x^k_{i_k}-y^k_{i_k})$ to replace $\gamma_k(x^k_{i_k}-x^k_{i_k})$ in scheme \eqref{schemes} and then updates $y^k_{i_k}=x^k_{i_k}$ with keeping other coordinates of $y^k$. In this scheme, the inertial term can be active for all iterations. However, the convergence of such algorithm is beyond the proof techniques proposed in this paper, and of course,  deserves further study.

\section{Conclusion}
In this paper, we studied the non-ergodic computational complexity of the Heavy-ball methods in the convex setting. Under different assumptions, we proved the non-ergodic sublinear and linear convergence rates for the algorithm, respectively. In both cases, we made much more relaxed assumptions than appeared in the existing literatures. Our proof was motivated by the analysis on a novel dynamical system. We extended our results to the multi-block coordinate descent Heavy-ball algorithm for both cyclic and stochastic update rules. The application to decentralized optimization demonstrated the advantage of our analysis techniques.\\

\textbf{Acknowledgments:} The authors are indebted to anonymous referees for their useful suggestions.
 We are grateful for the support from the the Major State Research Development Program (2016YFB0201305), and National Key Research and Development Program of China (2017YFB0202003), and National Natural Science Foundation of Hunan (2018JJ3616).

\emph{\onecolumn
\appendix
\begin{center}
{\large \bf Supplementary materials for \textit{Non-ergodic Convergence Analysis for Heavy-Ball Algorithm}}
\end{center}}

\section*{Proof of Lemma \ref{le1}}
By the scheme for updating $x^{k+1}$,
\begin{equation}\label{le1-t1}
    \frac{x^k-x^{k+1}}{\gamma_k}+\frac{\beta_k}{\gamma_k}(x^k-x^{k-1})=\nabla f(x^k).
\end{equation}
By Lipschitz continuity of $\nabla f$,
\begin{align}\label{le1-t3}
    f(x^{k+1})-f(x^k)&\leq \langle \nabla f(x^k),x^{k+1}-x^k\rangle+\frac{L}{2}\|x^{k+1}-x^k\|^2.
\end{align}
Combining (\ref{le1-t1}) and (\ref{le1-t3}), we have
\begin{align}\label{le1-t4}
    f(x^{k+1})-f(x^{k})&\overset{(\ref{le1-t1})+(\ref{le1-t3})}{\leq} \frac{\beta_k}{\gamma_k}\langle x^{k}-x^{k-1}, x^{k+1}-x^k\rangle+\left(\frac{L}{2}-\frac{1}{\gamma_k}\right)\|x^{k+1}-x^k\|^2\nonumber\\
    &\overset{a)}{\leq}\frac{\beta_k}{2\gamma_k}\|x^{k}-x^{k-1}\|^2+\left(\frac{L}{2}-\frac{1}{\gamma_k}+\frac{\beta_k}{2\gamma_k}\right)\|x^{k+1}-x^k\|^2.
 \end{align}
where $a)$ uses the Cauchy-Schwarz inequality $\langle x^{k}-x^{k-1}, x^{k+1}-x^k\rangle\leq \frac{1}{2}\|x^{k}-x^{k-1}\|^2+\frac{1}{2}\|x^{k+1}-x^k\|^2$. A simple calculation gives
\begin{align}\label{le1-t5}
     &\left[f(x^{k})+\frac{\beta_k}{2\gamma_k}\|x^{k}-x^{k-1}\|^2\right]-\left[f(x^{k+1})+\frac{\beta_{k}}{2\gamma_{k}}\|x^{k+1}-x^k\|^2\right]\geq\left(\frac{1-\beta_k}{\gamma_k}-\frac{L}{2}\right)\|x^{k+1}-x^k\|^2.
\end{align}
Since $(\beta_k)_{k\geq 0}$ is non-increasing, so is $\left(\frac{\beta_k}{2\gamma_k}=\frac{\beta_k L}{4(1-\beta_k)c}\right)_{k\geq 0}$,
\begin{align}\label{le1-t6}
f(x^{k+1})+\frac{\beta_{k}}{2\gamma_{k}}\|x^{k+1}-x^k\|^2
  \geq f(x^{k+1})+\frac{\beta_{k+1}}{2\gamma_{k+1}}\|x^{k+1}-x^k\|^2.
\end{align}
Summing (\ref{le1-t5}) and (\ref{le1-t6}), we obtain the
(\ref{nresult}).

\section*{Proof of Lemma \ref{lem-sub}}
By a direct computation and Lemma \ref{le1}, we have
\begin{align}\label{sketch-1}
  \xi_k- \xi_{k+1}\geq\frac{L(1-c)}{4c}\cdot(\|x^{k+1}-x^{k}\|^2+\|x^{k}-x^{k-1}\|^2).
\end{align}
The convexity of $f$ implies that
\begin{equation}\label{lem-sub-t5}
    f(x^{k})-f(\overline{x^{k}})\leq \langle\nabla f(x^{k}),x^{k}-\overline{x^{k}}\rangle.
\end{equation}
Summing (\ref{le1-t1}) and (\ref{lem-sub-t5})  yields
\begin{align}\label{lem-sub-t6}
    f(x^{k})-f(\overline{x^{k}})
    &\leq \frac{\beta_k}{\gamma_k}\langle x^k-x^{k-1},x^{k}-\overline{x^{k}}\rangle+ \left\langle  \frac{x^k-x^{k+1}}{\gamma_k},x^{k}-\overline{x^{k}}\right\rangle\nonumber\\
    &\overset{a)}{\leq}\frac{\beta_k}{\gamma_k}\|x^k-x^{k-1}\|\cdot\|x^{k}-\overline{x^{k}}\|+\frac{1}{\gamma_k}\|x^{k+1}-x^k\|\cdot\|x^{k}-\overline{x^{k}}\|\nonumber\\
    &\overset{b)}{\leq}\frac{1}{\gamma_k}\left(\|x^{k+1}-x^k\|+\|x^k-x^{k-1}\|\right)\times\|x^{k}-\overline{x^{k}}\|,
\end{align}
where $a)$ is due to the Cauchy-Schwarz inequality, $b)$ is due to the fact $0\leq\beta_k<1$. Combining (\ref{Lyapunov}) and (\ref{lem-sub-t6}), we have
\begin{align}\label{lem-sub-t7}
    \xi_{k}&\leq\frac{1}{\gamma_k}\left(\|x^{k+1}-x^k\|+\|x^k-x^{k-1}\|\right)\cdot\|x^{k}-\overline{x^{k}}\|+\delta_{k}\|x^k-x^{k-1}\|^2.\nonumber
\end{align}
Let
 \begin{align*}
    a^k&:=\left(
                              \begin{array}{c}
                              \frac{1}{\gamma_k}\|x^{k+1}-x^k\| \\
                              \frac{1}{\gamma_k}\|x^k-x^{k-1}\| \\
                               \delta_{k}\|x^k-x^{k-1}\| \\
                              \end{array}
                            \right),\quad
                            b^k:=\left(
                              \begin{array}{c}
                              \|x^{k}-\overline{x^{k}}\| \\
                              \|x^{k}-\overline{x^{k}}\| \\
                               \|x^k-x^{k-1}\| \\
                              \end{array}
                            \right).
\end{align*}
Using this and the definition of $\xi_{k}$ (\ref{Lyapunov}), we have:
\begin{align}
    (\xi_{k})^2\leq \left|\langle a^k,b^k\rangle\right|^2\leq\|a^k\|^2\cdot\|b^k\|^2.
\end{align}
A direct calculation give
\begin{align*}
    \|a^k\|^2\leq(\delta_{k}^2+\frac{1}{\gamma^2_k})\times(\|x^{k+1}-x^k\|^2+\|x^k-x^{k-1}\|^2)
\end{align*}
and
\begin{equation*}
    \|b^k\|^2\leq 2\|x^{k}-\overline{x^{k}}\|^2+\|x^k-x^{k-1}\|^2.
\end{equation*}
Thus we have
\begin{align}\label{sketch-2}
   \xi_{k}^2\leq(\delta_{k}^2+\frac{1}{\gamma^2_k})\cdot(\|x^{k+1}-x^k\|^2+\|x^k-x^{k-1}\|^2)\times(2\|x^{k}-\overline{x^{k}}\|^2+\|x^k-x^{k-1}\|^2).
\end{align}
Combining (\ref{sketch-1}) and (\ref{sketch-2}) completes the proof.

\section*{Proof of Theorem \ref{sub-n}}
 By Lemma \ref{lem-sub},
 $$\sup_{k}\{\xi_{k}\}<+\infty.$$
We also have
$$\sup_{k}\{f(x^k)\}<+\infty$$
 and
 $$\sup_{k}\{\|x^{k+1}-x^k\|^2\}<+\infty.$$
By the coercivity of $F$, sequences $(x^k)_{k\geq 0}$ and $(\overline{x^k})_{k\geq 0}$ are bounded. So $R<+\infty$. By the assumptions on $\gamma_k$ and $\beta_k$,
  $$\sup_k\{\varepsilon_k\}<+\infty.$$
  It is easy to see $\|x^{k}-x^{k-1}\|^2\leq 2R$, then,
\begin{equation*}
    \sup_{k}\{\varepsilon_k(2\|x^{k}-\overline{x^{k}}\|^2+\|x^{k}-x^{k-1}\|^2)\}\leq 4R\cdot \sup_k\{\varepsilon_k\}.
\end{equation*}
By Lemma \ref{lem-sub},  we then have
\begin{equation*}
    \xi_{k}^2\leq 4R\cdot \sup_k\{\varepsilon_k\}(\xi_k-\xi_{k+1}).
\end{equation*}
Since $\xi_{k+1}\leq\xi_{k}$, we have
\begin{equation*}
    \xi_{k+1}\xi_{k}\leq 4R\cdot \sup_k\{\varepsilon_k\}(\xi_k-\xi_{k+1}).
\end{equation*}
That is also
\begin{equation*}
    \frac{1}{\xi_{i+1}}-\frac{1}{\xi_{i}}\geq \frac{1}{4R\cdot \sup_k\{\varepsilon_k\}}.
\end{equation*}
Summing the inequality from $i=0$ to $k$ gives
\begin{equation*}
    \xi_k\leq \frac{1}{\frac{1}{\xi(0)}+\frac{k}{4R\cdot \sup_k\{\varepsilon_k\}}}.
\end{equation*}
Since $f(x^k)-\min f\leq \xi_k$, we then complete the proof.

\section*{Proof of Corollary \ref{lem-b}}
We just need to verify the boundedness of the points.
Let $x^*$ be a minimizer of $f$.
Noting $I-\gamma_k\nabla f(\cdot)$ is  contractive when $0<\gamma_k<\frac{2}{L}$,
\begin{align*}
\|x^{k+1}-x^*\|&= \|[x^k-\gamma_k \nabla f(x^k)+\beta_k(x^k-x^{k-1})]-[x^*-\gamma_k \nabla f(x^*)]\|\nonumber\\
&\quad\leq\|[x^k-\gamma_k \nabla f(x^k)]-[x^*-\gamma_k \nabla f(x^*)]\|\nonumber\\
&\quad+\|\beta_k(x^k-x^*+x^*-x^{k-1})\|\nonumber\\
&\quad\leq (1+\beta_k)\|x^k-x^*\|+\beta_k\|x^{k-1}-x^*\|.
\end{align*}
Denote that
 $$t_k:=\|x^k-x^*\|,\quad h_k:=t_k+\beta_{k-1} t_{k-1}.$$
 We then have
 \begin{equation}\label{app2-assump}
    t_{k+1}\leq (1+\beta_k)t_k+\beta_k t_{k-1}.
\end{equation}
Adding $\beta_k t_k$ to both sides of (\ref{app2-assump}),
 \begin{align}\label{app2-t1}
    t_{k+1}+\beta_k t_k&\leq (1+\beta_k)t_k+\beta_k t_{k-1}+\beta_k t_k\nonumber\\
    &\leq (1+2\beta_k)(t_k+\beta_k t_{k-1}).
 \end{align}
 Noting the decent of $(\beta_k)_{k\geq 0}$, (\ref{app2-t1}) is actually
  \begin{equation*}
    t_{k+1}+\beta_k t_k\leq (1+2\beta_k)(t_k+\beta_{k-1} t_{k-1}).
 \end{equation*}
We then have
   \begin{equation*}
    h_{k+1}\leq (1+2\beta_k)h_k\leq e^{2\beta_k}h_k.
 \end{equation*}
 Thus, for any $k$
 \begin{equation*}
    h_{k+1}\leq  e^{2\sum_{i=1}^k\beta_i}h_1<+\infty.
 \end{equation*}
 The boundedness of $\{h_k\}_{k\geq 0}$ directly yields the boundedness of $\{t_k\}_{k\geq 0}$.

\section*{Proof of Theorem \ref{thm-linear}}
With (\ref{condtion}), we have
\begin{equation*}
    2\|x^{k}-\overline{x^{k}}\|^2\leq \frac{2}{\nu}( f(x^{k})-\min f)\leq \frac{2}{\nu}\xi_k.
\end{equation*}
On the other hand, from the definition of (\ref{Lyapunov}),
\begin{equation*}
    \|x^{k}-x^{k-1}\|^2\leq  \frac{1}{\delta_k}  \xi_{k}.
\end{equation*}
With Lemma \ref{lem-sub}, we then derive
\begin{equation*}
     \xi_{k}^2\leq \varepsilon_k\left(\frac{1}{\delta_k}+\frac{2}{\nu}\right)(\xi_k-\xi_{k+1})\cdot\xi_k.
\end{equation*}
Thus, we define
$$\ell:=\sup_{k}\left\{\varepsilon_k\left(\frac{1}{\delta_k}+\frac{2}{\nu}\right)\right\}<+\infty.$$
And then, we have the following result,
\begin{equation*}
    \xi_{k}\xi_{k+1}\leq \ell(\xi_k-\xi_{k+1})\cdot\xi_k
\end{equation*}
due to that $0\leq \xi_{k+1}\leq \xi_{k}$.
With basic algebraic computation,
\begin{equation*}
    \frac{ \xi_{k+1}}{ \xi_{k}}\leq \frac{\ell}{1+\ell}.
\end{equation*}
By defining $\omega=\frac{\ell}{1+\ell}$, we then prove the result.

\section*{Proof of Lemma \ref{md-le1}}
For any $i\in[1,2,\ldots,m]$,
\begin{equation}\label{md-le1-t1}
    \frac{x^k_i-x^{k+1}_i}{\gamma_{k,i}}+\frac{\beta_{k,i}}{\gamma_{k,i}}(x^k_i-x^{k-1}_i)=\nabla_i^k f.
\end{equation}
With (\ref{lpm}), we can have
\begin{align}\label{md-le1-t3}
    &f(x_1^{k+1},\ldots,x_{i-1}^{k+1},x_i^{k+1},x_{i+1}^{k}\ldots,x_m^k)-f(x_1^{k+1},\ldots,x_{i-1}^{k+1},x_i^{k},x_{i+1}^{k}\ldots,x_m^k)\nonumber\\
    &\quad\quad\leq \langle \nabla_i^k f,x^{k+1}_i-x^k_i\rangle+\frac{L}{2}\|x^{k+1}_i-x^k_i\|^2.
\end{align}
Combining (\ref{md-le1-t1}) and (\ref{md-le1-t3}),
\begin{align}\label{md-le1-t4}
    &f(x_1^{k+1},\ldots,x_{i-1}^{k+1},x_i^{k+1},x_{i+1}^{k}\ldots,x_m^k)-f(x_1^{k+1},\ldots,x_{i-1}^{k+1},x_i^{k},x_{i+1}^{k}\ldots,x_m^k)\nonumber\\
    &\quad\overset{(\ref{md-le1-t1})+(\ref{md-le1-t3})}{\leq} \frac{\beta_{k,i}}{\gamma_{k,i}}\langle x^{k}_i-x^{k-1}_i, x^{k+1}_i-x^k_i\rangle+\left(\frac{L}{2}-\frac{1}{\gamma_{k,i}}\right)\|x^{k+1}_i-x^k_i\|^2\nonumber\\
    &\quad\overset{a)}{\leq}\frac{\beta_{k,i}}{2\gamma_{k,i}}\|x^{k}_i-x^{k-1}_i\|^2+\left(\frac{L_i}{2}-\frac{1}{\gamma_{k,i}}+\frac{\beta_{k,i}}{2\gamma_{k,i}}\right)\|x^{k+1}_i-x^k_i\|^2.
 \end{align}
where $a)$ uses the Schwarz inequality $\langle x^{k}_i-x^{k-1}_i, x^{k+1}_i-x^k_i\rangle\leq \frac{1}{2}\|x^{k}_i-x^{k-1}_i\|^2+\frac{1}{2}\|x^{k+1}_i-x^k_i\|^2$.
Summing (\ref{md-le1-t4}) from $i=1$ to $m$,
\begin{align}\label{md-le1-t5}
f(x^{k+1})-f(x^k)\leq\sum_{i=1}^m\frac{\beta_{k,i}}{2\gamma_{k,i}}\|x^{k}_i-x^{k-1}_i\|^2+\sum_{i=1}^m\left(\frac{L_i}{2}-\frac{1}{\gamma_{k,i}}+\frac{\beta_{k,i}}{2\gamma_{k,i}}\right)\|x^{k+1}_i-x^k_i\|^2.
 \end{align}
With direct calculations and the non-increasity of $(\beta_{k,i})_{k\geq 0}$, we then obtain (\ref{nresult}).

\section*{Proof of Lemma \ref{lem-sub-m}}
With  Lemma \ref{md-le1}, direct computing yields
\begin{align}\label{sketch-1-m}
    \hat{\xi}_k- \hat{\xi}_{k+1}&\geq\sum_{i=1}^m \frac{1}{2}(\frac{1-\beta_{k,i}}{\gamma_{k,i}}-\frac{L_i}{2})\cdot(\|x^{k+1}_i-x^k_i\|^2+\|x^{k}_i-x^{k-1}_i\|^2)\nonumber\\
    &=\frac{\underline{L}}{4}(\frac{1}{c}-1)\cdot(\|x^{k+1}-x^k\|^2+\|x^{k}-x^{k-1}\|^2).
\end{align}
With the convexity of $H$, we then have
\begin{align}\label{lem-sub-t6-m}
    &f(x^{k})-f(\overline{x^{k}})\leq \langle\nabla H(x^{k}),x^{k}-\overline{x^{k}}\rangle\nonumber\\
    &\quad=\sum_{i=1}^m\frac{\beta_{k,i}}{\gamma_{k,i}}\langle x^{k}_i-x^{k-1}_i,x^{k}_i-\left[\overline{x^{k}}\right]_i\rangle+\sum_{i=1}^m\langle  \frac{x^k_i-x^{k+1}_i}{\gamma_{k,i}},x^{k}_i-\left[\overline{x^{k}}\right]_i\rangle\nonumber\\
    &\quad+\sum_{i=1}^m\langle \nabla_i f(x^{k})-\nabla_i^k f,x^{k}_i-\left[\overline{x^{k}}\right]_i\rangle\nonumber\\
    &\quad\overset{a)}{\leq}\sum_{i=1}^m\frac{\beta_{k,i}}{\gamma_{k,i}}\|x^{k}_i-x^{k-1}_i\|\cdot\|x^{k}_i-\left[\overline{x^{k}}\right]_i\|
    +\sum_{i=1}^m \frac{1}{\gamma_{k,i}}\|x^{k+1}_i-x^k_i\|\cdot\|x^{k}_i-\left[\overline{x^{k}}\right]_i\|\nonumber\\
    &\quad+\sum_{i=1}^m L\|x^{k+1}-x^k\|\cdot\|x^{k}_i-\left[\overline{x^{k}}\right]_i\|\nonumber\\
    &\quad\overset{b)}{\leq}\sum_{i=1}^m\left(\frac{\|x^{k+1}_i-x^k_i\|}{\gamma_{k,i}}+\frac{\|x^{k}_i-x^{k-1}_i\|}{\gamma_{k,i}}+L\|x^{k+1}-x^k\|\right)\cdot\|x^{k}_i-\left[\overline{x^{k}}\right]_i\|,
\end{align}
where $a)$ is due to the Schwarz inequalities and the smooth assumption \textbf{A1}, $b)$ depends on the fact $0\leq\beta_{k,i}<1$. With (\ref{Lyapunov}) and (\ref{lem-sub-t6}), we have
\begin{align}\label{lem-sub-t7-m}
    \hat{\xi}_{k}&\leq\sum_{i=1}^m\left(\frac{\|x^{k+1}_i-x^k_i\|}{\gamma_{k,i}}+\frac{\|x^{k}_i-x^{k-1}_i\|}{\gamma_{k,i}}+L\|x^{k+1}-x^k\|\right)\cdot\|x^{k}_i-\left[\overline{x^{k}}\right]_i\|\nonumber\\
    &+\sum_{i=1}^m \delta_{k,i}\|x^{k}_i-x^{k-1}_i\|^2.\nonumber
\end{align}
Let
 \begin{align*}
    \hat{a}^k:=\left(
                              \begin{array}{c}
                              \frac{1}{\gamma_{k,1}}\|x^{k+1}_1-x^k_1\| \\
                              \vdots\\
                              \frac{1}{\gamma_{k,m}}\|x^{k+1}_m-x^k_m\| \\
                              \frac{1}{\gamma_{k,1}}\|x^{k}_1-x^{k-1}_1\| \\
                              \vdots \\
                              \frac{1}{\gamma_{k,m}}\|x^{k}_m-x^{k-1}_m\| \\
                              L\|x^{k+1}-x^k\|\\
                              \vdots\\
                              L\|x^{k+1}-x^k\|\\
                               \delta_{k,1}\|x^{k}_1-x^{k-1}_1\| \\
                               \vdots \\
                                \delta_{k,m}\|x^{k}_m-x^{k-1}_m\| \\
                              \end{array}
                            \right),\quad \hat{b}^k:=\left(
                              \begin{array}{c}
                              \|x^{k}_1-\left[\overline{x^{k}}\right]_1\| \\
                              \vdots\\
                               \|x^{k}_m-\left[\overline{x^{k}}\right]_m\| \\
                                 \|x^{k}_1-\left[\overline{x^{k}}\right]_1\| \\
                              \vdots\\
                               \|x^{k}_m-\left[\overline{x^{k}}\right]_m\| \\
                                 \|x^{k}_1-\left[\overline{x^{k}}\right]_1\| \\
                              \vdots\\
                               \|x^{k}_m-\left[\overline{x^{k}}\right]_m\| \\
                               \|x^{k}_1-x^{k-1}_1\| \\
                               \vdots\\
                               \|x^{k}_m-x^{k-1}_m\| \\
                              \end{array}
                            \right).
\end{align*}
Using this and the definition of $\hat{\xi}_{k}$ (\ref{Lyapunov-m}), we have:
\begin{align*}
    (\hat{\xi}_{k})^2\leq\left|\langle \hat{a}^k,\hat{b}^k\rangle\right|^2\leq\|\hat{a}^k\|^2\cdot\|\hat{b}^k\|^2.
\end{align*}
Direct calculation yields
\begin{equation*}
    \|\hat{a}^k\|^2\leq\max\{\sum_{i=1}^m\left(\delta_{k,i}^2+\frac{1}{\gamma^2_{k,i}}\right),m\cdot L^2\}\cdot(\|x^{k+1}-x^k\|^2+\|x^{k}-x^{k-1}\|^2)
\end{equation*}
and
\begin{equation*}
    \|\hat{b}^k\|^2\leq 2\|x^{k+1}-\overline{x^{k+1}}\|^2+\|x^{k+1}-x^k\|^2.
\end{equation*}
Thus, we derive
\begin{align}\label{sketch-2-m}
    (\hat{\xi}_{k})^2&\leq\max\{\sum_{i=1}^m\left(\delta_{k,i}^2+\frac{1}{\gamma^2_{k,i}}\right),m\cdot L^2\}\nonumber\\
    &\times(\|x^{k+1}-x^k\|^2+\|x^{k}-x^{k-1}\|^2)\times(2\|x^{k+1}-\overline{x^{k+1}}\|^2+\|x^{k+1}-x^k\|^2).
\end{align}
Combining (\ref{sketch-1-m}) and (\ref{sketch-2-m}), we then prove the result.
\section*{Proof of Theorem \ref{them-sub-m}}
With Lemma \ref{lem-sub-m},  $\sup_{k}\{\hat{\xi}_{k}\}<+\infty$, thus, $\sup_{k}\{f(x^k)\}<+\infty$ and $\sup_{k}\{\|x^{k}-x^{k-1}\|^2\}<+\infty$. Noting the coercivity of $f$, sequences $(x^k)_{k\geq 0}$ and $(\overline{x^k})_{k\geq 0}$ are  bounded. With the assumptions on $\gamma_{k,i}$ and $\beta_{k,i}$, $\sup_k\{\hat{\varepsilon}_k\}<+\infty$. Noting that $\|x^{k-1}-x^{k}\|^2\leq 2R$,
\begin{equation*}
    \hat{\varepsilon}_k(2\|x^{k+1}-\overline{x^{k+1}}\|^2+\|x^{k-1}-x^{k}\|^2)\leq 4R\cdot \hat{\varepsilon}_k.
\end{equation*}
With Lemma \ref{lem-sub-m},  we then have
\begin{equation*}
    \hat{\xi}_{k}\hat{\xi}_{k+1}\leq\hat{\xi}_{k}^2\leq 4R\cdot \hat{\varepsilon}_k(\hat{\xi}_k-\hat{\xi}_{k+1}).
\end{equation*}
Then, we have
\begin{equation*}
    \hat{\xi}_k\leq \frac{4R\cdot \sup_{k}\{\hat{\varepsilon}_k\}}{k}.
\end{equation*}
Using the fact $f(x^k)-\min f\leq \hat{\xi}_k$, we then obtain the result.
\section*{Proof of Theorem \ref{them-linear-m}}

With the restricted strong convexity, we have
\begin{equation*}
    2\|x^{k}-\overline{x^{k}}\|^2\leq \frac{2}{\nu}\hat{\xi}_{k}.
\end{equation*}
The direct computing yields
\begin{equation*}
    \|x^{k-1}-x^{k}\|^2\leq\frac{\hat{\xi}_{k}}{\min_i\{\delta_{k,i}\}}
\end{equation*}
With Lemma \ref{lem-sub-m},
\begin{equation*}
    (\hat{\xi}_{k})^2\leq\left(\hat{\varepsilon}_k+\frac{2}{\nu}+\frac{1}{\min_i\{\delta_{k,i}\}}\right)(\hat{\xi}_k-\hat{\xi}_{k+1})\hat{\xi}_{k}.
\end{equation*}
It is easy to see that $\ell:=\sup_{k}\big\{\hat{\varepsilon}_k+\frac{2}{\nu}+\frac{1}{\min_i\{\delta_{k,i}\}}\big\}$. Then, we have
\begin{equation*}
    f(x^k)-\min f\leq \hat{\xi}_k= O\left(\left(\frac{\ell}{1+\ell}\right)^k\right).
\end{equation*}
Letting $\omega=\frac{\ell}{1+\ell}$, we then prove the result.

\section*{Proof of Lemma \ref{ms-le1}}
In the following, the sub-algebra $\chi^k$ is defined as
\begin{equation}\label{sub-al}
    \chi^k:=\sigma(x^0,x^1,\ldots,x^k).
\end{equation}
In the $k$-th iteration,
\begin{equation}\label{ms-le1-t1}
    \frac{x^k_{i_k}-x^{k+1}_{i_k}}{\gamma_{k}}+\frac{\beta_{k}}{\gamma_{k}}(x^k_{i_k}-x^{k-1}_{i_k})=\nabla_{i_k} f(x^k).
\end{equation}
With the Lipschitz of $\nabla f$, we can have
\begin{eqnarray}\label{ms-le1-t3}
    f(x^{k+1})-f(x^k)\leq \langle -\nabla_{i_k} f(x^k),x^k_{i_k}-x^{k+1}_{i_k}\rangle+\frac{L}{2}\|x^{k+1}-x^k\|^2,
\end{eqnarray}
where we used the fact $\langle -\nabla_{i_k} f(x^k),x^k_{i_k}-x^{k+1}_{i_k}\rangle=\langle -\nabla f(x^k),x^k-x^{k+1}\rangle$.
Combining (\ref{ms-le1-t1}) and (\ref{ms-le1-t3}),
\begin{eqnarray}\label{ms-le1-t4}
    &&f(x^{k+1})-f(x^k)\overset{(\ref{ms-le1-t1})+(\ref{ms-le1-t3})}{\leq} \frac{\beta_{k}}{\gamma_{k}}\langle x^{k}_{i_k}-x^{k-1}_{i_k}, x^{k+1}_{i_k}-x^k_{i_k}\rangle+(\frac{L}{2}-\frac{1}{\gamma_{k}})\|x^{k+1}-x^k\|^2\nonumber\\
    &&~\quad\quad\overset{a)}{\leq}\frac{\sqrt{m}\beta_{k}}{2\gamma_{k}}\|x^{k}_{i_k}-x^{k-1}_{i_k}\|^2+\left(\frac{L}{2}-\frac{1}{\gamma_{k}}+\frac{\beta_{k}}{2\sqrt{m}\gamma_{k}}\right)\|x^{k+1}-x^k\|^2.
 \end{eqnarray}
where $a)$ uses the Schwarz inequality $\langle x^{k}_i-x^{k-1}_i, x^{k+1}_i-x^k_i\rangle\leq \frac{\sqrt{m}}{2}\|x^{k}_i-x^{k-1}_i\|^2+\frac{1}{2\sqrt{m}}\|x^{k+1}_i-x^k_i\|^2$ and the fact $\|x^{k+1}-x^k\|^2=\|x^{k+1}_{i_k}-x^k_{i_k}\|^2$.
Taking conditional conditional expectations of (\ref{ms-le1-t4}) on $\chi^k$,
\begin{align}\label{ms-le1-t5}
&\EE[f(x^{k+1})\mid\chi^k]-f(x^k)\leq\frac{\beta_{k}}{2\sqrt{m}\gamma_{k}}\|x^{k}-x^{k-1}\|^2+(\frac{L}{2}-\frac{1}{\gamma_{k}}+\frac{\beta_{k}}{2\sqrt{m}\gamma_{k}})\EE(\|x^{k+1}-x^k\|^2\mid\chi^k).
 \end{align}
Taking total expectations on (\ref{ms-le1-t5}), and using $\EE(\EE(\cdot\mid\chi^k))=\EE(\cdot)$,
\begin{eqnarray*}\label{ms-le1-t6}
\EE f(x^{k+1}) -\EE f(x^k)\leq\frac{\beta_{k}}{2\sqrt{m}\gamma_{k}}\EE\|x^{k}-x^{k-1}\|^2+\left(\frac{L}{2}-\frac{1}{\gamma_{k}}+\frac{\beta_{k}}{2\sqrt{m}\gamma_{k}}\right)\EE\|x^{k+1}-x^k\|^2.
 \end{eqnarray*}
 Thus, we have
 \begin{align*}\label{ms-le1-t7}
    & \left[\EE f(x^{k})+\frac{\beta_k}{2\sqrt{m}\gamma_k}\EE\|x^{k}-x^{k-1}\|^2\right]-\left[\EE f(x^{k+1})+\frac{\beta_{k}}{2\sqrt{m}\gamma_{k}}\EE\|x^{k+1}-x^{k}\|^2\right]\\
     &\quad\quad\geq\left(\frac{1-\beta_k/\sqrt{m}}{\gamma_k}-\frac{L}{2}\right)\EE\|x^{k+1}-x^k\|^2.\nonumber
\end{align*}
With the non-increasity of $(\beta_k)_{k\geq 0}$, $\left(\frac{\beta_k}{2\sqrt{m}\gamma_k}=\frac{\beta_k L}{4(1-\beta_k/\sqrt{m})c}\right)_{k\geq 0}$ is also non-increasing, and then we prove the result.
\section*{Proof of Lemma \ref{lem-sub-s}}
With  Lemma \ref{ms-le1},
\begin{eqnarray}\label{sketch-1-ms}
 \EE\bar{\xi}_k-\EE\bar{\xi}_{k+1}\geq\frac{L}{4}(\frac{1}{c}-1)\cdot(\EE\|x^{k+1}-x^k\|^2+\EE\|x^{k}-x^{k-1}\|^2).
\end{eqnarray}
The convexity of $H$ yields
\begin{equation}\label{lem-sub-ms-t5}
    f(x^{k})-f(\overline{x^{k}})\leq \langle\nabla f(x^{k}),x^{k}-\overline{x^{k}}\rangle\leq \|\nabla f(x^{k})\|\cdot\|x^{k}-\overline{x^{k}}\|.
\end{equation}
With (\ref{Lyapunov-s}), we have
\begin{eqnarray}\label{lem-sub-ms-t7}
    \bar{\xi}_{k}\leq \|\nabla f(x^{k})\|\cdot\|x^{k}-\overline{x^{k}}\|+\delta_{k}\|x^k-x^{k-1}\|^2.
\end{eqnarray}
Denote that
 \begin{align*}
    \bar{a}^k:=\left(
                              \begin{array}{c}
                              \|\nabla f(x^{k})\| \\
                               \delta_{k}\|x^k-x^{k-1}\| \\
                              \end{array}
                            \right),
                            \bar{b}^k:=\left(
                              \begin{array}{c}
                              \|x^{k}-\overline{x^{k}}\| \\
                               \|x^k-x^{k-1}\| \\
                              \end{array}
                            \right).
\end{align*}
Thus, we have
\begin{equation}\label{lem-sub-ms-t1}
    (\EE\bar{\xi}_{k})^2\leq (\EE\langle \bar{a}^k,\bar{b}^k\rangle)^2\leq(\EE\|\bar{a}^k\|\cdot\|\bar{b}^k\|)^2\leq \EE\|\bar{a}^k\|^2\cdot\EE\|\bar{b}^k\|^2.
\end{equation}
With the scheme of the algorithm,
\begin{align}\label{lem-sub-t6-}
    &\EE\|\nabla_{i_k} f(x^k)\|^2=\EE\|\frac{x^k_{i_k}-x^{k+1}_{i_k}}{\gamma_{k}}+\frac{\beta_{k}}{\gamma_{k}}(x^k_{i_k}-x^{k-1}_{i_k})\|^2\nonumber\\
    &\quad\leq \frac{2\EE\|x^{k+1}-x^k\|^2+2\EE\|x^{k}-x^{k-1}\|^2}{\gamma_{k}^2}
\end{align}
where we used the fact $0\leq\beta_k<1$.
Direct calculation yields
\begin{align*}
    \EE\|\bar{a}^k\|^2&=\EE\|\nabla f(x^{k}) \|^2+\delta_k^2\cdot\EE\|x^k-x^{k-1}\|^2\nonumber\\
    &\leq m\cdot\EE\|\nabla_{i_k} f(x^{k}) \|^2+\delta_k^2\cdot\EE\|x^k-x^{k-1}\|^2\nonumber\\
    &\leq(\delta_{k}^2+\frac{2m}{\gamma^2_k})\times(\EE\|x^{k+1}-x^k\|^2+\EE\|x^k-x^{k-1}\|^2)
\end{align*}
and
\begin{equation*}
    \EE\|\bar{b}^k\|^2\leq \EE\|x^{k}-\overline{x^{k}}\|^2+\EE\|x^k-x^{k-1}\|^2.
\end{equation*}
Therefore, we derive
\begin{align}\label{sketch-2-ms}
    (\EE\bar{\xi}_{k})^2&\leq(\delta_{k}^2+\frac{2m}{\gamma^2_k})\times(\EE\|x^{k+1}-x^k\|^2+\EE\|x^k-x^{k-1}\|^2)\nonumber\\
    &\times(\|x^{k}-\overline{x^{k}}\|^2+\|x^k-x^{k-1}\|^2).
\end{align}
Combining (\ref{sketch-1-ms}) and (\ref{sketch-2-ms}), we then prove the result.
\section*{Proof of Theorem \ref{ms-th1}}
With Lemma \ref{ms-le1}, we can see
\begin{equation}\label{ms-th1-t1}
    \sum_{k}\EE\|x^{k+1}-x^k\|^2<+\infty.
\end{equation}
Noting $i_k$ is selected with uniform probability   from $\{1,2,\ldots,m\}$,
\begin{equation}\label{ms-th1-t2}
    \EE(\|\nabla_{i_k} f(x^k)\|^2\mid\chi^k)=\frac{\|\nabla f(x^k)\|^2}{m}.
\end{equation}
Taking total expectations on both sides of (\ref{ms-th1-t2}) and using $\EE(\EE(\cdot\mid\chi^k))=\EE(\cdot)$,
\begin{equation}\label{ms-th1-t2.5}
    \EE(\|\nabla_{i_k} f(x^k)\|^2)=\frac{\EE\|\nabla f(x^k)\|^2}{m}.
\end{equation}
On the other hand, the scheme of the algorithm gives
\begin{align}\label{ms-th1-t3}
    &\EE(\|\nabla_{i_k} f(x^k)\|^2\mid\chi^k)=\|\frac{x^k_{i_k}-x^{k+1}_{i_k}}{\gamma_{k}}+\frac{\beta_{k}}{\gamma_{k}}(x^k_{i_k}-x^{k-1}_{i_k})\|^2\nonumber\\
    &\qquad\leq\frac{2\|x^{k+1}-x^k\|^2}{\gamma_k}+\frac{2\|x^{k}-x^{k-1}\|^2}{\gamma_k}\nonumber\\
    &\qquad\leq\frac{2\|x^{k+1}-x^k\|^2}{\gamma_0}+\frac{2\|x^{k}-x^{k-1}\|^2}{\gamma_0}.
\end{align}
We also take total expectations on both sides of (\ref{ms-th1-t3}),
\begin{align}\label{ms-th1-t4}
    \EE\|\nabla_{i_k} f(x^k)\|^2\leq\frac{2\EE\|x^{k+1}-x^k\|^2}{\gamma_0}+\frac{2\EE\|x^{k}-x^{k-1}\|^2}{\gamma_0}.
\end{align}
Combining (\ref{ms-th1-t4}) and (\ref{ms-th1-t2.5}),
\begin{align}\label{ms-th1-t5}
    \EE\|\nabla f(x^k)\|^2\leq\frac{2m\EE\|x^{k+1}-x^k\|^2}{\gamma_0}+\frac{2m\EE\|x^{k}-x^{k-1}\|^2}{\gamma_0}.
\end{align}
With (\ref{ms-th1-t5}) and (\ref{ms-th1-t1}),
\begin{align}\label{ms-th1-t6}
    \sum_{k}\EE\|\nabla f(x^k)\|^2<+\infty.
\end{align}
Then, we have
\begin{equation*}
     \min_{0\leq i\leq k}\EE\|\nabla f(x^k)\|^2=o\left(\frac{1}{k}\right).
\end{equation*}
Finally, using the  Schwarz inequality, we then prove the result.

\section*{Proof of Theorem \ref{ms-th1-f}}
It is easy to see
$$\bar{\varepsilon}_k\cdot(\EE\|x^{k}-\overline{x^{k}}\|^2+\EE\|x^{k}-x^{k-1}\|^2)\leq 4R\cdot\sup_{k}\{\bar{\varepsilon}_k\}$$
 With Lemma \ref{ms-th1-t6}, we have
\begin{equation}
    (\EE\bar{\xi}_{k})^2\leq 4R\cdot\sup_{k}\{\bar{\varepsilon}_k\}(\EE\bar{\xi}_k-\EE\bar{\xi}_{k+1}).
\end{equation}
Considering the fact $0\leq \EE\bar{\xi}_{k+1}\leq \EE\bar{\xi}_{k}$, we have
\begin{equation}
    \EE\bar{\xi}_{k} \EE\bar{\xi}_{k+1}\leq 4R\cdot\sup_{k}\{\bar{\varepsilon}_k\}(\EE\bar{\xi}_k-\EE\bar{\xi}_{k+1}).
\end{equation}
Then we have
\begin{equation}
    \EE\bar{\xi}_{k}=O\left(\frac{4R\cdot\sup_{k}\{\bar{\varepsilon}_k\}}{k}\right).
\end{equation}
Considering that $\EE f(x^k)-\min f\leq  \EE\bar{\xi}_{k}$, we then prove the result.
\section*{Proof of Theorem \ref{th:linear-convergence}}
With the restricted strong convexity, we have
\begin{equation*}
    \EE\|x^{k+1}-\overline{x^{k+1}}\|^2\leq\frac{1}{\nu}\EE\bar{\xi}_{k}.
\end{equation*}
The direct computing yields
\begin{equation*}
    \EE\|x^{k-1}-x^{k}\|^2\leq\frac{\EE\bar{\xi}_{k}}{\bar{\delta}_{k}}
\end{equation*}
With Lemma \ref{ms-th1-f},
\begin{equation*}
    (\EE\bar{\xi}_{k})^2\leq\left(\bar{\varepsilon}_k+\frac{1}{\nu}+\frac{1}{\bar{\delta}_{k}}\right)\cdot(\EE\bar{\xi}_k-\EE\bar{\xi}_{k+1})\cdot\EE\bar{\xi}_{k}.
\end{equation*}
It is easy to see that $\bar{\varepsilon}_k+\frac{1}{\nu}+\frac{1}{\bar{\delta}_{k}}$ is bounded by some positive constant $\bar{\ell}>0$. With the fact $\EE\bar{\xi}_{k+1}\leq \EE\bar{\xi}_{k}$,
\begin{equation*}
    \EE\bar{\xi}_{k+1}\cdot\EE\bar{\xi}_{k}\leq\bar{\ell}\cdot(\EE\bar{\xi}_k-\EE\bar{\xi}_{k+1})\cdot\EE\bar{\xi}_{k}.
\end{equation*}
 Then, we have
\begin{equation*}
    \EE f(x^k)-\min f\leq \EE\hat{\xi}_k= O\left(\left(\frac{\bar{\ell}}{1+\bar{\ell}}\right)^k\right).
\end{equation*}
Letting $\omega=\frac{\bar{\ell}}{1+\bar{\ell}}$, we then prove the result.

\end{document}